\documentclass[11pt]{article}
\usepackage[T1]{fontenc}
\usepackage[utf8]{inputenc}
\usepackage{lmodern}
\usepackage{amsmath,amssymb,amsthm,mathrsfs}
\usepackage{mathtools}
\usepackage[margin=1in]{geometry}
\usepackage{microtype}
\usepackage{parskip}
\usepackage{needspace}
\usepackage{xcolor}

\definecolor{revisionpurple}{RGB}{128,0,128}

\usepackage{enumitem}
\usepackage[hidelinks,bookmarksnumbered=true]{hyperref}

\hypersetup{pdftitle={Packing sets by rigid motions and similarities},pdfauthor={Daewoong Cheong, Thang Pham, and Chun-Yen Shen}}
\numberwithin{equation}{section}
\newtheorem{theorem}{Theorem}[section]
\newtheorem{lemma}[theorem]{Lemma}
\newtheorem{proposition}[theorem]{Proposition}
\newtheorem{corollary}[theorem]{Corollary}
\theoremstyle{remark}
\newtheorem{remark}[theorem]{Remark}
\newtheorem{example}[theorem]{Example}
\newcommand{\R}{\mathbb R}
\newcommand{\Z}{\mathbb Z}
\newcommand{\dimH}{\dim_H}
\newcommand{\dimSF}{\dim_{\mathcal{SF}}}
\newcommand{\dimS}{\dim_{\mathcal S}}
\newcommand{\Leb}{\mathcal L}
\newcommand{\TT}{\mathbb T}
\newcommand{\wh}{\widehat}
\DeclareMathOperator{\spt}{spt}
\DeclareMathOperator{\dist}{dist}
\DeclareMathOperator{\Isom}{Isom}
\DeclareMathOperator{\Sim}{Sim}
\allowdisplaybreaks[1]
\title{Packing sets by rigid motions and similarities}
\author{Daewoong Cheong \and Thang Pham \and Chun-Yen Shen}
\date{}
\begin{document}
\maketitle
\begin{abstract}
For nonempty Borel sets \(E\subset\mathbb R^d\) and \(\Theta\subset O(d)\times\mathbb R^d\), where \(d\ge2\), let \(\Theta(E)\) be the union of the rigid copies \(gE+z\), where \((g,z)\in\Theta\). We prove that \(\Theta(E)\) has positive Lebesgue measure whenever
\[
\dimSF E+\dim_H\Theta>\frac{d(d+1)}2,
\]
where \(\dimSF E\) is the \emph{spherical Fourier dimension} of \(E\). This resolves a conjecture on packing by arbitrary Borel families of rigid motions. Our proof develops a discrete-to-continuous framework
by adapting a finite-field averaging argument to \(\mathbb R^d\)
using Sobolev estimates on the group of rigid motions. We also obtain
analogous packing results for similarities. The same framework can be
used to show that if \(\tau\) is a compactly supported probability
measure on the similarity group with finite
\(\bigl(\frac{d(d-1)}{2}+1\bigr)\)-energy, then every compactly
supported \(\tau\)-stationary probability measure on \(\mathbb R^d\)
has a smooth density. We deduce that a nonempty compact set has
nonempty interior if its family of similarity self-embeddings has
Hausdorff dimension greater than \(\frac{d(d-1)}{2}+1\).
\end{abstract}
\noindent\textbf{2020 Mathematics Subject Classification.} 28A75.

\noindent\textbf{Keywords.} Packing sets, rigid motions, similarity transformations, Hausdorff dimension, spherical averages, incidences.

\tableofcontents
\section{Introduction}\label{sec:introduction}

\begingroup
Let $d\ge2$, and let $O(d)$ denote the orthogonal group. We identify the space $\Isom(\R^d)$ of rigid motions with $O(d)\times\R^d$. For $\theta=(g,z)$, set $\theta(x)=gx+z$. Given Borel sets $E\subset\R^d$ and $\Theta\subset\Isom(\R^d)$, define
\[
\Theta(E)=\bigcup_{\theta\in\Theta}\theta(E)
=\{gx+z:x\in E,\ (g,z)\in\Theta\}.
\]
The \emph{packing problem} asks when a set of Lebesgue measure zero can contain a copy $\theta(E)$ for every $\theta\in\Theta$. We seek conditions ensuring that $\Leb^d(\Theta(E))>0$, and we also obtain lower bounds for $\dimH\Theta(E)$. Throughout, we write $\Leb^d$ for \emph{Lebesgue measure}, $\mathcal H^u$ for \emph{$u$-dimensional Hausdorff measure}, and $\dimH$ for \emph{Hausdorff dimension}. \emph{Haar probability measure} on $O(d)$ is denoted by $dg$, and $B(x,r)$ is the closed ball with centre $x$ and radius $r$. Set
\[
m=\frac{d(d-1)}2,\qquad n=m+d=\frac{d(d+1)}2.
\]
These are the dimensions of $O(d)$ and $\Isom(\R^d)$, respectively. On $\Isom(\R^d)$, we use the product of a fixed bi-invariant Riemannian metric on $O(d)$ and the Euclidean metric on $\R^d$. Other parameter manifolds carry fixed smooth Riemannian metrics. Hausdorff dimension does not depend on the choice of smooth Riemannian metric.

For a subset $A$ in a Euclidean space or a smooth manifold, $\mathcal M(A)$ denotes the set of all nonzero finite positive Radon measures with compact support contained in $A$. We write $\spt\lambda$ for the support of a measure $\lambda$. We normalize these measures to have mass one when convenient. For a measurable map $F$, the \emph{pushforward measure} $F_\#\mu$ is defined by $F_\#\mu(B)=\mu(F^{-1}(B))$ for Borel sets $B$. We write $X\lesssim Y$ if $X\le CY$ for a constant independent of the scale under consideration, and $X\sim Y$ if both inequalities hold. Dependence on fixed parameters is indicated when needed. A finite positive measure $\lambda$ is called \emph{$s$-Frostman} if $\lambda(B(x,r))\lesssim r^s$ for every $x$ and $0<r\le1$.
\endgroup

For a finite measure $\mu$ on $\R^d$, its \emph{Fourier transform} is
\[
\wh\mu(\xi)=\int e^{-2\pi i x\cdot\xi}\,d\mu(x),\qquad \xi\in\R^d.
\]
Let $\sigma_{d-1}$ be the normalized surface measure on $S^{d-1}$. We define the \emph{spherical average}
\begin{equation}\label{eq:spherical-average}
\sigma(\mu)(r)=\int_{S^{d-1}}|\wh\mu(r\omega)|^2\,d\sigma_{d-1}(\omega).
\end{equation}
Mattila \cite{Mattila87} studied spherical averages in connection with distance sets and intersections. Following \cite{IMPPSS}, we define the \emph{spherical Fourier dimension} of a nonempty set $E$ by
\[
\dimSF E=\sup\bigl\{\beta\in[0,d]:\exists\mu\in\mathcal M(E),\ 
\sigma(\mu)(r)\lesssim(1+r)^{-\beta}\text{ for }r\ge0\bigr\}.
\]
We will frequently use the following decay condition
for a compactly supported measure $\mu$ and some $0\le\beta\le d$:
\begin{equation}\label{eq:spherical}
\sigma(\mu)(r)\lesssim(1+r)^{-\beta},\qquad r\ge0.
\end{equation}

The \emph{Sobolev dimension} of a nonempty set $E\subset\R^d$ is defined by
\begin{equation}\label{eq:sobolev-dimension}
\dimS E=\sup\left\{u\in\R:\exists\mu\in\mathcal M(E),\ 
\int_{\R^d}|\wh\mu(\xi)|^2(1+|\xi|)^{u-d}\,d\xi<\infty\right\}.
\end{equation}
For $a\in\R$, we write $H^a$ for the $L^2$-based
\emph{Sobolev space of order $a$}. Its norms and the relations between energy, Sobolev regularity, and Hausdorff dimension are recalled in Section~\ref{sec:preliminaries}.

Packing problems include classical questions about unions of lines, circles, and planes; see Mattila \cite{MattilaBook} and Keleti \cite{Keleti}. Besicovitch constructed a planar set of Lebesgue measure zero containing a unit segment in every direction. For circles, Wolff \cite{WolffCircles} proved that the union has positive area if the set of centers has Hausdorff dimension greater than one, even when the radii vary. Oberlin \cite{OberlinPlanes} obtained packing results for arbitrary families of affine planes. Ham, Ko, Lee, and Oh \cite{HKLO} studied translates and dilates of nondegenerate curves using local smoothing estimates.

\begingroup
Iosevich, Mattila, Palsson, Pham, Pham, Senger, and Shen \cite{IMPPSS} considered general Borel sets. Their Theorem~1.4 treats product families $\Theta=G\times Z$, under the additional assumption $\dimH G>(d-1)(d-2)/2$. They proved that
\[
\dimSF E+\dimH G+\dimH Z>\frac{d(d+1)}2
\quad\Longrightarrow\quad
\Leb^d(GE+Z)>0,
\]
and obtained corresponding dimension bounds and results for product families of similarities. Their Conjecture~1 asks whether the rigid-motion conclusion holds for arbitrary Borel $\Theta\subset\Isom(\R^d)$, with $\dimH\Theta$ replacing $\dimH G+\dimH Z$. Our Theorem~\ref{thm:main} proves this conjecture. For translations and dilations, \cite[Theorems 1.2 and 1.3]{IMPPSS} already allow arbitrary Borel parameter families under suitable dimension hypotheses.
\endgroup

Over finite fields, the second author \cite{Pham} proved packing results for arbitrary families of rigid motions. In particular, if $q$ is a prime power with $q\equiv3\pmod4$, then \cite[Theorem 1.11]{Pham} gives
\begin{equation}\label{eq:finite}
|\Theta(E)|\gtrsim\min\left\{q^2,\frac{|E|^{1/2}|\Theta|}{q}\right\}
\end{equation}
for $E\subset\mathbb F_q^2$ and $\Theta\subset O(2,\mathbb F_q)\times\mathbb F_q^2$. The authors of \cite{IMPPSS} provided another proof and noted that a planar Euclidean analogue remained to be obtained. The planar case of Corollary~\ref{cor:hausdorff} gives the corresponding positive-measure conclusion.

\subsection{Main results}
We first state the rigid-motion theorem.

\begin{theorem}\label{thm:main}
Let $E\subset\R^d$ and $\Theta\subset \Isom(\R^d)$ be nonempty Borel sets. Then the following hold.
\begin{enumerate}
\item[(i)] One has
\[
\dimS\Theta(E)\ge\dimSF E+\dimH\Theta-\frac{d(d-1)}2.
\]
\item[(ii)] If $\dimSF E+\dimH\Theta>d(d+1)/2$, then $\Leb^d(\Theta(E))>0$.
\item[(iii)] For $0<u<d$, if $\dimSF E+\dimH\Theta>d(d-1)/2+u$, then $\dimH\Theta(E)\ge u$.
\end{enumerate}
\end{theorem}

For any $\theta\in\Theta$, the inclusion $\theta(E)\subset\Theta(E)$ gives $\dimH\Theta(E)\ge\dimH\theta(E)=\dimH E$. Combining this elementary bound with Theorem~\ref{thm:main} gives
\[
\dimH\Theta(E)\ge\max\left\{\dimH E,
\min\left\{d,\dimSF E+\dimH\Theta-\frac{d(d-1)}2\right\}\right\}.
\]
For every fixed nonempty finite set, the preceding Hausdorff bound is attained at every possible dimension of the motion family, even when $E$ affinely spans $\R^d$; see Proposition~\ref{prop:finite-source}.
The example in Proposition~\ref{prop:sharp} shows that the strict inequality in part~(ii) cannot be replaced by equality.

Note that the standard spherical-average estimates imply
\begin{equation}\label{eq:sf-lower}
\dimSF E\ge\frac{d-1}{d}\dimH E;
\end{equation}
see \cite[Lemma 1.1]{IMPPSS}. Thus, Theorem~\ref{thm:main} gives the following consequence depending only on Hausdorff dimension.

\Needspace{6\baselineskip}
\begin{corollary}\label{cor:hausdorff}
Let $E\subset\R^d$ and $\Theta\subset \Isom(\R^d)$ be nonempty Borel sets.
\begin{enumerate}
\item[(i)] If $\frac{d-1}{d}\dimH E+\dimH\Theta>d(d+1)/2$, then $\Leb^d(\Theta(E))>0$.
\item[(ii)] For $0<u<d$, if $\frac{d-1}{d}\dimH E+\dimH\Theta>d(d-1)/2+u$, then $\dimH\Theta(E)\ge u$.
\end{enumerate}
\end{corollary}

For $d=2$, part~(i) reads
\begin{equation}\label{eq:planar-target}
\frac12\dimH E+\dimH\Theta>3
\quad\Longrightarrow\quad\Leb^2(\Theta(E))>0.
\end{equation}
The coefficient $1/2$ cannot be increased in a statement valid for all $E$ and $\Theta$; see Proposition~\ref{prop:coefficient}. When $0 < \dimH E < 1$, the classical circular-average estimates give a stronger sufficient condition. For $0\le s\le2$, define
\begin{equation}\label{eq:planar-beta}
\beta_2(s)=
\begin{cases}
s,&0\le s\le1/2,\\
1/2,&1/2\le s\le1,\\
s/2,&1\le s\le2.
\end{cases}
\end{equation}
The estimates in \cite[Section 2.1]{IMPPSS}, including Wolff's circular-average estimate \cite{Wolff}, give $\dimSF E\ge\beta_2(\dimH E)$. Hence, $\beta_2(\dimH E)+\dimH\Theta>3$ implies positive area.

As another planar application, let $K$ be a compact planar \emph{Kakeya set}, that is, a compact subset of $\R^2$ containing a unit line segment in every direction. Oberlin's Fourier-decay result \cite[Proposition 2 and the following remark]{OberlinKakeya}, also recovered in \cite[Corollary 2.4]{FHK}, implies that for every $0<\beta<2$ there is a probability measure $\mu$ supported on $K$ such that $|\wh\mu(\xi)|\lesssim(1+|\xi|)^{-\beta/2}$. Averaging this bound over circles gives $\dimSF K=2$. Consequently, Theorem~\ref{thm:main}(ii) yields
\[
\dimH\Theta>1
\quad\Longrightarrow\quad
\Leb^2\!\left(\bigcup_{(g,z)\in\Theta}(gK+z)\right)>0
\]
for every Borel family $\Theta\subset \Isom(\R^2)$, even when $K$ itself has area zero.

We next consider \emph{similarity transformations}. We identify $\Sim(\R^d)$ with $O(d)\times\R^d\times\R_{>0}$. For a Borel family $\Omega\subset\Sim(\R^d)$, define
\[
\Omega(E)=\{rgx+z:x\in E,\ (g,z,r)\in\Omega\}.
\]
Averaging over both rotations and scales allows the spherical Fourier dimension in Theorem~\ref{thm:main} to be replaced by Hausdorff dimension.

\Needspace{9\baselineskip}
\begin{theorem}\label{thm:similarity}
Let $E\subset\R^d$ and $\Omega\subset\Sim(\R^d)$ be nonempty Borel sets.
\begin{enumerate}
\item[(i)] One has
\[
\dimS\Omega(E)\ge\dimH E+\dimH\Omega-\frac{d^2-d+2}{2}.
\]
\item[(ii)] If $\dimH E+\dimH\Omega>(d^2+d+2)/2$, then $\Leb^d(\Omega(E))>0$.
\item[(iii)] For $0<u<d$, if $\dimH E+\dimH\Omega>(d^2-d+2)/2+u$, then $\dimH\Omega(E)\ge u$.
\end{enumerate}
\end{theorem}

Theorem~\ref{thm:similarity} removes the product assumption in \cite[Theorem 1.5]{IMPPSS}. Example~\ref{ex:similarity} shows that equality does not suffice in part~(ii).

We now return to rigid motions in the plane. When $E$ is contained in a line, the geometry of $E$ can
give stronger bounds. Suppose that $E\subset\R^2$ and
$\Theta\subset\Isom(\R^2)$ are nonempty Borel sets, $E$ lies in an affine
line, $s=\dimH E\in(0,1]$, and $t=\dimH\Theta>1$.
Proposition~\ref{prop:line}, obtained from the Ren--Wang theorem
\cite{RenWang}, gives
\[
 \dimH\Theta(E)\ge
 \min\left\{s+t-1,\frac{3s+t-1}{2},s+1\right\}.
\]
If $\mathcal H^1(E)>0$, Proposition~\ref{prop:line-area} gives
\[
 \dimH\Theta>2\quad\Longrightarrow\quad\Leb^2(\Theta(E))>0.
\]
The strict threshold $2$ is sharp; see Example~\ref{ex:line-area}.
This improves the sufficient condition $\dimH\Theta>5/2$ supplied by
Theorem~\ref{thm:main} and the general bound $\beta_2(1)=1/2$.

We also obtain several incidence estimates as applications of our
main results. In Theorem~\ref{thm:incidence}, we establish a
quantitative incidence bound for rigid motions in every dimension
$d\ge2$. The examples in Section~\ref{subsec:incidence-sharpness}
show that the bound is sharp in all three regimes, including the
logarithmic factor at the critical exponent. In addition,
Proposition~\ref{prop:point-line} recovers Orponen's point-line
incidence inequality from our Sobolev estimate. 

Besides the packing and incidence results, our averaging
estimates also have applications to stationary measures.
We now turn to these applications.

\subsection{Stationary measures and similarity self-embeddings}

\begingroup
Recall that
\[
 m=\dim O(d)=\frac{d(d-1)}2,\qquad
 n=\dim\Isom(\R^d)=m+d=\frac{d(d+1)}2.
\]
\endgroup

Let \(\Phi\subset\operatorname{Sim}(\mathbb R^d)\) be a compact family of similarities, and let \(\mu\) be a compactly supported probability measure such that
\[
\Phi(\spt\mu)\subset \spt\mu.
\]
We are interested in the regularity of \(\mu\) when it is preserved by averaging over this family. More precisely, given a probability measure \(\tau\) on \(\Phi\), we call \(\mu\) \textit{stationary} if
\[
P_\tau\mu:=\int_\Phi\theta_\#\mu\,d\tau(\theta)=\mu.
\]
Our next result gives Frostman and energy conditions on \(\tau\) under which every compactly supported stationary probability measure has a \(C^\infty\) density. As a consequence, we obtain a criterion for a compact set admitting many similarity self-embeddings to have nonempty interior.

\begin{theorem}\label{app-thm-stationary}
Let $d\ge1$, and let $\tau$ be a compactly supported probability measure
on $\Sim(\R^d)$.  If $\tau$ is $t$-Frostman with $t>m+1$, then every
compactly supported stationary probability measure for $P_\tau$ has a
$C^\infty$ density.  The conclusion also holds if the Frostman hypothesis
is replaced by
\begin{equation}\label{app-eq-critical-law}
 \tau\in H^{-d/2}(\Sim(\R^d)).
\end{equation}
\end{theorem}

Since $\dim\Sim(\R^d)=n+1$, condition
\eqref{app-eq-critical-law} is equivalent to
the finiteness of the $(m+1)$-energy of $\tau$.
More precisely, let $d_{\Sim}$ be a distance induced
locally by a smooth Riemannian metric on each component
of $\Sim(\R^d)$, with the two components separated
by a positive distance. Then this condition is equivalent to
\[
\iint d_{\Sim}(\theta,\theta')^{-(m+1)}
\,d\tau(\theta)\,d\tau(\theta')<\infty.
\]

We now consider similarities that map a compact set
into itself. Let $E\subset\R^d$ be nonempty and compact.
A similarity $\theta$ is called a
\emph{similarity self-embedding} of $E$
if $\theta(E)\subset E$. We denote the family of
all such similarities by
\[
\Sigma_E
=
\{(g,z,r)\in\Sim(\R^d):rgE+z\subset E\}.
\]

\begin{corollary}\label{app-cor-self-embeddings}
Let $d\ge1$, and let $E\subset\R^d$ be nonempty and compact.  If a compact
family $\Phi\subset\Sigma_E$ satisfies $\dimH\Phi>m+1$, then $E$ has
nonempty interior.  The same conclusion holds if $\Phi$ supports a
probability measure with finite $(m+1)$-energy in local parameter
coordinates.  Consequently, if $E$ has empty interior, then
\begin{equation}\label{app-eq-self-embedding-bound}
 \dimH\Sigma_E\le m+1,
\end{equation}
and $\Sigma_E$ supports no compactly supported probability measure of
finite $(m+1)$-energy.
\end{corollary}

For $d\ge2$, the strict dimension hypothesis in Corollary~\ref{app-cor-self-embeddings} already implies $\Leb^d(E)>0$ by Theorem~\ref{thm:similarity}. Indeed, $\Phi(E)\subset E$ and part~(iii) give
\[
\dimH E\ge\min\{d,\dimH E+\dimH\Phi-(m+1)\}.
\]
Since $\dimH\Phi>m+1$, this inequality forces $\dimH E=d$.
Hence, $\dimH E+\dimH\Phi>n+1$, so part~(ii) applies. The stationary-measure argument strengthens this conclusion to nonempty interior and also covers the finite-energy condition.

Let $\Phi\subset\Sim(\R^d)$ be a nonempty compact family of
contracting similarities. The \emph{attractor} of $\Phi$ is the
unique nonempty compact set $E\subset\R^d$ satisfying
\[
E=\Phi(E)=\bigcup_{\theta\in\Phi}\theta(E).
\]
Indeed, since $\Phi$ is compact and every similarity in $\Phi$
has contraction ratio less than one, there is a constant
$0<\rho<1$ such that
\[
|\theta(x)-\theta(y)|\le \rho|x-y|
\]
for all $\theta\in\Phi$ and $x,y\in\R^d$.
Consequently, the map $A\mapsto\Phi(A)$ is a contraction
on the space of nonempty compact subsets of $\R^d$,
equipped with the Hausdorff metric. The existence and
uniqueness of $E$ then follow from the Banach fixed-point theorem.

Since every $\theta\in\Phi$ maps $E$ into itself,
Corollary~\ref{app-cor-self-embeddings} implies that
$E$ has nonempty interior whenever $\dimH\Phi>m+1$.
The dimension threshold is sharp without further assumptions
on $E$. The proofs and the sharpness example are given in
Section~\ref{sec:new-applications}.

Related questions concerning attractors and the regularity of stationary measures have been studied in several settings. For example, Hochman \cite[Section 1.1]{Hochman} studies whether attractors of compact, possibly uncountable families of similarities of the line have nonempty interior. Persson \cite[Theorem 1]{Persson} obtained smooth stationary densities for random affine contractions of the line with smooth additive noise. \'{S}wi\k{a}tkowski \cite[Sections 3 and 4]{Swiatkowski} proved absolute continuity under geometric or Fourier assumptions on convolution powers of the affine parameter law. Here the hypothesis concerns the dimension or energy of a possibly singular joint law of rotations, scales, and translations.

For isometries, the existence of a compactly supported stationary measure
forces a common fixed point. Indeed, suppose that
$\tau$ is a compactly supported probability measure on $\Isom(\R^d)$
and that $P_\tau\mu=\mu$ for a compactly supported probability measure
$\mu$. Let $\bar x=\int x\,d\mu(x)$, then stationarity gives
\begin{align*}
0
&=\iint\bigl(|gx+z-\bar x|^2-|x-\bar x|^2\bigr)
        \,d\mu(x)d\tau(g,z)\\
&=\int|g\bar x+z-\bar x|^2\,d\tau(g,z).
\end{align*}
The second equality uses $|g(x-\bar x)|=|x-\bar x|$ and
$\int(x-\bar x)\,d\mu(x)=0$, which eliminates the cross term.
The final integrand is nonnegative and continuous, so it vanishes on
$\spt\tau$. Hence, every motion in $\spt\tau$ fixes $\bar x$. Thus,
$\delta_{\bar x}$ is also stationary and has no density, while
\[
 \spt\tau\subset\{(g,\bar x-g\bar x):g\in O(d)\}.
\]
In particular, a $t$-Frostman isometry law with $t>m$ admits no
compactly supported stationary probability measure.

A similar restriction holds for isometries mapping a compact set
into itself. Indeed, if $E\subset\R^d$ is nonempty and compact,
then every isometry $\theta$ satisfying $\theta(E)\subset E$
must satisfy $\theta(E)=E$. To see this, suppose that
$x\in E\setminus\theta(E)$ and let
$\varepsilon=\operatorname{dist}(x,\theta(E))>0$.
For $0\le i<j$, we have
\[
|\theta^i x-\theta^j x|
=|x-\theta^{j-i}x|\ge\varepsilon.
\]
Since $\theta(E)\subset E$, all iterates $\theta^j x$ lie in
$E$. The inequality above makes them an infinite
$\varepsilon$-separated subset, contradicting compactness. Thus,
$\theta(E)=E$.
Such an isometry fixes the unique center of the smallest ball
containing $E$. Since the stabilizer in
$\Isom(\R^d)$ a point in $\mathbb{R}^d$ is isomorphic to $O(d)$, it follows that
\[
\dimH\{\theta\in\Isom(\R^d):\theta(E)\subset E\}\le m.
\]
This bound is sharp, as shown by $E=S^{d-1}$,
which has an empty interior. Related questions concerning sets preserved by families
of linear transformations over finite fields have also
been studied in~\cite{PHS}.

\subsection{Main ideas of the proofs:
from finite fields to Euclidean space}
\begingroup
We now return to the packing problem and explain
the main ideas of our proofs. Our starting point is the proof of Theorem~1.6 in \cite[Section 5]{IMPPSS}. For $E\subset\mathbb F_q^2$ and $\Theta\subset O(2,\mathbb F_q)\times\mathbb F_q^2$, with $q\equiv3\pmod4$, let
\[
\lambda_\Theta(y)=\sum_{x\in E}\sum_{(g,z)\in\Theta}
\mathbf1_{\{y=gx+z\}}.
\]
Cauchy--Schwarz gives
\[
|\Theta(E)|\ge\frac{|E|^2|\Theta|^2}{\sum_y\lambda_\Theta(y)^2}.
\]
Taking the Fourier transform in the translation variable and applying Cauchy--Schwarz separates a spherical Fourier average of $E$ from the partial Fourier transforms of $\mathbf1_\Theta$. A spherical restriction estimate controls the first factor. Plancherel controls the second by $|\Theta|$. These estimates imply \eqref{eq:finite}, since
\[
\sum_y\lambda_\Theta(y)^2
\lesssim q^{-2}|E|^2|\Theta|^2+q|E|^{3/2}|\Theta|.
\]
\endgroup

A direct extension to $\R^d$ encounters a difficulty at the last step. Let $\mu\in\mathcal M(E)$ satisfy \eqref{eq:spherical}, let $\tau\in\mathcal M(\Theta)$ be $t$-Frostman, and let $\nu$ be the image of $\mu\times\tau$ under $(x,g,z)\mapsto gx+z$. Then
\[
\wh\nu(\xi)=\int e^{-2\pi i z\cdot\xi}\wh\mu(g^{-1}\xi)\,d\tau(g,z).
\]
The measure $\tau$ may be singular with respect to $dg\,dz$. Thus, the partial Plancherel argument does not directly apply. For $\Theta=G\times Z$, the proof in \cite{IMPPSS} chooses separate Frostman measures on $G$ and $Z$. Their product separates the rotation and translation variables, allowing a rotational average to be combined with an energy estimate for the translation measure.

To overcome this difficulty, we study the averaging operator
on the full motion group, following the operator approach
suggested by the finite-field argument in Appendix~B.
Related Sobolev estimates for averaging operators have appeared
in packing problems for curves and hypersurfaces
\cite{HKLO,ILT}, while spherical Fourier averages have been
used to study convolution operators on the motion group
\cite{BGTT}. In our setting, the key step is to obtain
Sobolev estimates in both the rotation and translation
variables. More precisely, define
\[
T_\mu f(g,z)=\int f(gx+z)\,d\mu(x).
\]
Partial Plancherel and Haar averaging first give
\[
\|T_\mu f\|_{L^2(\Isom(\R^d))}\lesssim\|f\|_{H^{-\beta/2}(\R^d)}.
\]
We then show that smooth weights preserve spherical Fourier decay. Derivatives in $g$ insert polynomial weights in $x$ and powers of the spatial frequency. The weighted estimate therefore controls rotation derivatives, and interpolation yields
\[
\|T_\mu f\|_{H^a(\Isom(\R^d))}\lesssim_a
\|f\|_{H^{a-\beta/2}(\R^d)},\qquad a\ge0.
\]
This estimate includes derivatives in both $g$ and $z$. The joint Frostman condition implies $\tau\in H^{-\beta/2}(\Isom(\R^d))$ when $\beta+t>d(d+1)/2$. Thus, Sobolev duality gives
\[
\left|\int f\,d\nu\right|
=|\langle\tau,T_\mu f\rangle|
\lesssim\|\tau\|_{H^{-\beta/2}}\|f\|_2.
\]
Hence, $\nu$ has an $L^2$ density. Sobolev duality replaces the final Plancherel step and uses the dimension of the full motion family, without separate assumptions on its rotation and translation projections. Whenever $2a>n-t$, duality gives $\nu\in H^{\beta/2-a}(\R^d)$, which yields the Sobolev and Hausdorff dimension bounds.

The proof for similarities follows the same approach, but the
additional scale parameter allows us to obtain a stronger result.
Indeed, averaging over both rotations and scales reduces the
Fourier estimate to an integral over an annulus, which can be
controlled by the standard Fourier estimate for Frostman measures.
Thus, unlike the rigid-motion case, we only need the Frostman
condition on the measure supported on $E$, without assuming
spherical Fourier decay. This gives the corresponding Sobolev
estimate for the similarity averaging operator. We then apply
Sobolev duality with a measure on the full similarity parameter
space to prove Theorem~\ref{thm:similarity}.

The rest of the paper is organized as follows. Section~\ref{sec:preliminaries} collects the analytic
preliminaries needed in the proofs. In
Sections~\ref{sec:rigid-proof} and~\ref{sec:similarity-proof},
we prove the main packing results for rigid motions
and similarities. Section~\ref{sec:new-applications}
establishes the stationary-measure and self-embedding
results using an extension of the similarity estimate.
In Section~\ref{sec:examples}, we give examples showing
the sharpness of the packing bounds.
Section~\ref{sec:applications} is devoted to quantitative
incidence estimates and their applications, including
point-line incidences, measure zero sets, and patterns contained
in a line. Appendix~\ref{app:planar-fourier}
provides an alternative proof of the planar estimate
using angular Fourier series. Appendix \ref{app:finite-field} presents a finite-field proof that inspires our arguments in $\mathbb{R}^d$.

\section{Preliminaries}\label{sec:preliminaries}
In this section, we recall the relations between Frostman bounds, energy, and Sobolev regularity. We also state two Fourier estimates used in the packing proofs.

\subsection{Frostman measures, energies, and Fourier averages}
Frostman's lemma implies that, for a nonempty Borel set $A$ in a Euclidean space or a smooth manifold and $0<t<\dimH A$, there is a $t$-Frostman probability measure in $\mathcal M(A)$; see \cite[Theorem 2.7]{MattilaBook}. 
\begingroup
A nonempty compact set $A$ is called \emph{$t$-Ahlfors regular} if it supports a probability measure $\lambda$ such that
\[
\lambda(B(x,r))\sim r^t,
\qquad x\in A,\quad 0<r\le1.
\]
\endgroup

For $0<u<d$ and a finite positive measure $\lambda$ on $\R^d$, the \emph{$u$-energy} is
\begin{equation}\label{eq:riesz-energy}
I_u(\lambda)=\iint|x-y|^{-u}\,d\lambda(x)d\lambda(y)
=c_{d,u}\int|\wh\lambda(\xi)|^2|\xi|^{u-d}\,d\xi.
\end{equation}
If $I_u(\lambda)<\infty$ and $\lambda\in\mathcal M(A)$, then $\dimH A\ge u$; see \cite[Chapter 3]{MattilaBook}. A $t$-Frostman measure has finite $u$-energy for $0<u<t$.

Similarly, if $\mu$ satisfies \eqref{eq:spherical}, then polar coordinates give
{
\[
I_u(\mu)\lesssim\mu(\R^d)^2+\int_1^\infty r^{u-1}\sigma(\mu)(r)\,dr<\infty,
\qquad 0<u<\min\{\beta,d\}.
\]
}
Thus, $\dimSF E\le\dimH E$ for every nonempty Borel set $E\subset\R^d$.

We will use the following standard Fourier estimate on balls; see \cite[Section 3.8]{MattilaBook}.
\begin{lemma}\label{lem:ball}
Let $\mu$ be a compactly supported $s$-Frostman measure on $\R^d$, with $0\le s\le d$. For every bounded function $w$ on $\spt\mu$,
\[
\int_{|\xi|\le R}|\wh{w\mu}(\xi)|^2\,d\xi
\lesssim\|w\|_\infty^2 R^{d-s},\qquad R\ge1.
\]
\end{lemma}

For completeness, we include its short Gaussian proof in Section~\ref{sec:similarity-proof}. The constants may depend on the Frostman constant and the total mass. In particular, this estimate applies to the smooth weights arising from derivatives of $\wh\mu$.

\subsection{Sobolev spaces and parameter measures}
\begingroup
For $a\in\R$, the Sobolev space $H^a(\R^d)$ consists of the tempered distributions $f$ for which
\[
\|f\|_{H^a(\R^d)}^2
=\int_{\R^d}|\wh f(\xi)|^2(1+4\pi^2|\xi|^2)^a\,d\xi<\infty.
\]
If $\lambda\in\mathcal M(A)$ and $\lambda\in H^0=L^2$, then $\Leb^d(\spt\lambda)>0$. For $0<u<d$, boundedness of $\wh\lambda$ near the origin and \eqref{eq:riesz-energy} give
\[
\lambda\in H^{(u-d)/2}
\quad\Longleftrightarrow\quad I_u(\lambda)<\infty
\quad\Longrightarrow\quad \dimH A\ge u.
\]
More generally, for any $u\in\R$, a measure $\lambda\in\mathcal M(A)\cap H^{(u-d)/2}$ is admissible in \eqref{eq:sobolev-dimension} at exponent $u$.
\endgroup

For the fixed bi-invariant Riemannian metric on $O(d)$, let $-\Delta_O$ be the nonnegative \emph{Laplace--Beltrami operator}. On $\Isom(\R^d)$, use the Sobolev norm defined by
\[
\mathcal A=1-\Delta_O-\Delta_z,
\qquad \|F\|_{H^a(\Isom(\R^d))}=\|\mathcal A^{a/2}F\|_2.
\]
These norms agree locally, up to equivalence, with the usual norms in smooth coordinate systemss. On $\Sim(\R^d)$, we only use functions and distributions supported in a compact subinterval of the positive scale axis. Extend them by zero to $O(d)\times\R^d\times\R$ and use the Sobolev norm defined by $1-\Delta_O-\Delta_z-\partial_r^2$. On other smooth manifolds, Sobolev regularity of compactly supported distributions is defined using coordinate systemss and a smooth partition of unity.

\begin{lemma}\label{lem:negative}
Let $M$ be a smooth $N$-dimensional manifold, and let $\tau$ be a compactly supported $t$-Frostman measure on $M$, where $0\le t\le N$. If $a>0$ and $2a>N-t$, then $\tau\in H^{-a}(M)$.
\end{lemma}

This is the negative-order Sobolev estimate needed to integrate the averaging operator against a measure on the parameter set. We prove it directly using the Bessel kernel. A stronger trace estimate is recorded in \cite[Theorem 2.1]{Dabrowski}.

\subsection{Smooth weights and spherical Fourier decay}
The next lemma controls the polynomial factors introduced by differentiation in the rotation variables.

\begin{lemma}\label{lem:weights}
Suppose that $\mu$ is a compactly supported finite measure satisfying \eqref{eq:spherical}. For every $w\in C_c^\infty(\R^d)$,
\[
\int_{ S^{d-1}}|\wh{w\mu}(r\omega)|^2\,d\sigma_{d-1}(\omega)
\le C_w(1+r)^{-\beta}.
\]
In particular, the same estimate holds for every partial derivative $\partial^\gamma\wh\mu$.
\end{lemma}

\section{Proof of Theorem~\ref{thm:main} and Corollary~\ref{cor:hausdorff}}\label{sec:rigid-proof}
In this section, we prove Theorem~\ref{thm:main} and Corollary~\ref{cor:hausdorff}. We first estimate the averaging operator in all motion variables, using spherical Fourier decay to control rotation derivatives. We then pair this estimate with a Frostman measure on the parameter family.

\subsection{The basic Fourier estimate}
Let $\mu$ be a compactly supported finite positive measure satisfying \eqref{eq:spherical}.
For a Schwartz function $f$, we first take the Fourier transform
of $T_\mu f$ only in $z$. Fubini's theorem and the change of variables
$y=gx+z$ give
\begin{equation}\label{eq:partial-transform}
\begin{aligned}
\mathcal F_z(T_\mu f)(g,\xi)
&=\wh f(\xi)\int e^{2\pi i gx\cdot\xi}\,d\mu(x)\\
&=\wh f(\xi)\wh\mu(-g^{-1}\xi),
\end{aligned}
\end{equation}
since $g^T=g^{-1}$. The minus sign has no effect on spherical averages.
Let $F_\xi(g)=\wh\mu(-g^{-1}\xi)$ and $\mathopen{\langle}{\xi}\mathclose{\rangle}=(1+4\pi^2|\xi|^2)^{1/2}$.
For $\xi\ne0$, the image of Haar probability measure under
$g\mapsto-g^{-1}\xi/|\xi|$ is normalized surface measure.
At $\xi=0$, both averages equal $\mu(\R^d)^2$. Thus,
\begin{equation}\label{eq:zeroth}
\|F_\xi\|_{L^2( O(d))}^2=\sigma(\mu)(|\xi|)\le C \mathopen{\langle}{\xi}\mathclose{\rangle}^{-\beta}.
\end{equation}
Thus, Plancherel's theorem gives
\begin{equation}\label{eq:basic-operator}
\|T_\mu f\|_{L^2(\Isom(\R^d))}^2
=\int_{\R^d}|\wh f(\xi)|^2\|F_\xi\|_2^2\,d\xi
\lesssim\|f\|_{H^{-\beta/2}(\R^d)}^2.
\end{equation}
The norm on the left uses the smooth measure $dg\,dz$. To integrate $T_\mu f$ against an arbitrary Frostman measure on $\Theta$, we need regularity in the joint variables $(g,z)$. We obtain it by controlling derivatives in the rotation variables.

\subsection{Rotation derivatives and the Sobolev estimate}
When we differentiate with respect to the rotation variables,
polynomial factors in $x$ appear.
We first show that these factors do not affect
the spherical Fourier decay of $\mu$.

\begin{proof}[Proof of Lemma~\ref{lem:weights}]
We use the convolution identity $\wh{w\mu}=\wh w*\wh\mu$ to
transfer the spherical decay of $\mu$ to $w\mu$. Choose $D>d+\beta$;
this ensures that the radial kernel remaining after the decay estimate
is integrable. Cauchy--Schwarz and the rapid decay of $\wh w$ give
\begin{align*}
|\wh{w\mu}(\xi)|^2
&=|\wh w*\wh\mu(\xi)|^2\\
&\le \|\wh w\|_1\int|\wh w(\xi-{\zeta})|\,|\wh\mu({\zeta})|^2\,d{\zeta}\\
&\le C_{D}\int(1+|\xi-{\zeta}|)^{-{D}}|\wh\mu({\zeta})|^2\,d{\zeta}.
\end{align*}
Let $k_D(v)=(1+|v|)^{-D}$. To recover the spherical average
of $\mu$, average the preceding inequality at $g\xi$ over $g\in O(d)$
and make the change of variables $\zeta\mapsto g\zeta$. Since $k_D$
is radial, this gives
\begin{align*}
\int_{O(d)}|\wh{w\mu}(g\xi)|^2\,dg
&\le C_{D}\int_{\R^d} k_{D}(\xi-\zeta)
\left(\int_{O(d)}|\wh\mu(g\zeta)|^2\,dg\right)d\zeta\\
&=C_{D}\int_{\R^d}k_{D}(\xi-\zeta)\sigma(\mu)(|\zeta|)\,d\zeta\\
&\le C_{D}\int_{\R^d}(1+|\xi-\zeta|)^{-{D}}(1+|\zeta|)^{-\beta}\,d\zeta.
\end{align*}
Since $1+|\xi|\le(1+|\xi-\zeta|)(1+|\zeta|)$, the last expression is at most
\[
C_{D}(1+|\xi|)^{-\beta}\int_{\R^d}(1+|v|)^{-{D}+\beta}\,dv
\le C(1+|\xi|)^{-\beta}
\]
\begingroup
by the choice $D>d+\beta$ and Haar averaging at $|\xi|=r$ gives the spherical average. For a multi-index $\gamma$, choose $w\in C_c^\infty(\R^d)$ equal to $(-2\pi i x)^\gamma$ near $\spt\mu$. Then
\[
\partial^\gamma\wh\mu=\wh{w\mu},
\]
so the same bound holds for every partial derivative. This completes the proof.
\endgroup
\end{proof}

The weighted estimate allows us to strengthen \eqref{eq:basic-operator} as follows.
\begin{proposition}\label{prop:operator}
Let $\mu$ be a compactly supported finite positive measure satisfying \eqref{eq:spherical}. Then, for every $a\ge0$ and every Schwartz function $f$ on $\R^d$,
\begin{equation}\label{eq:operator-estimate}
\|T_\mu f\|_{H^a(\Isom(\R^d))}
\le C_a\|f\|_{H^{a-\beta/2}(\R^d)}.
\end{equation}
In particular,
\[
T_\mu:L^2(\R^d)\longrightarrow H^{\beta/2}(\Isom(\R^d))
\]
is bounded.
\end{proposition}

\begin{proof}
Let $F_\xi$ and $\mathopen{\langle}{\xi}\mathclose{\rangle}$
be as in the preceding subsection, and set
$A_\xi=\mathopen{\langle}{\xi}\mathclose{\rangle}^{2}-\Delta_O$.
The partial Fourier transform in $z$ turns $\mathcal A$ into $A_\xi$.
Consequently, \eqref{eq:partial-transform} and Plancherel give
\[
\|T_\mu f\|_{H^a}^2
=\int_{\R^d}|\wh f(\xi)|^2
\|A_\xi^{a/2}F_\xi\|_2^2\,d\xi.
\]
It is therefore enough to prove
$\|A_\xi^{a/2}F_\xi\|_2^2
\lesssim_a\mathopen{\langle}{\xi}\mathclose{\rangle}^{2a-\beta}$
uniformly in $\xi$. We first estimate integer rotation derivatives,
then use the spectral theorem to obtain every real $a\ge0$.

\begingroup
To see how rotation derivatives enter, let $X$ be a skew-symmetric matrix. Differentiating along $g\exp(vX)$ gives
\[
\left.\frac{d}{dv}F_\xi(g\exp(vX))\right|_{v=0}
=\nabla\wh\mu(-g^{-1}\xi)\cdot Xg^{-1}\xi.
\]
Repeated differentiation gives a finite sum of derivatives $\partial^\gamma\wh\mu(-g^{-1}\xi)$ multiplied by polynomials in $g^{-1}\xi$. For a derivative of order $j$, both $|\gamma|$ and the polynomial degree are at most $j$.
The polynomial coefficients are bounded by
$C_j\mathopen{\langle}{\xi}\mathclose{\rangle}^{j}$. Apply
Lemma~\ref{lem:weights} to the finitely many derivatives
$\partial^\gamma\wh\mu$ and use Haar invariance. The equivalence of
integer Sobolev norms with norms of derivatives on $O(d)$ then gives
\endgroup
\begin{equation}\label{eq:derivatives}
\|F_\xi\|_{H^j( O(d))}^2\le C_j \mathopen{\langle}{\xi}\mathclose{\rangle}^{2j-\beta},\qquad j=0,1,2,\ldots.
\end{equation}

For $a=0$, the desired bound is \eqref{eq:zeroth}. Suppose
that $a>0$ and choose an integer $k>a$. Since the scalar
$\mathopen{\langle}{\xi}\mathclose{\rangle}^{2}$ commutes with
$-\Delta_O$, the spectral theorem and \eqref{eq:derivatives} yield
\[
\|A_{\xi}^{k/2}F_\xi\|_2^2
=\sum_{j=0}^k\binom{k}{j}\mathopen{\langle}{\xi}\mathclose{\rangle}^{2(k-j)}
\|(-\Delta_O)^{j/2}F_\xi\|_2^2
\le C_k\mathopen{\langle}{\xi}\mathclose{\rangle}^{2k-\beta}.
\]
Apply H\"older's inequality to the spectral measure of
$-\Delta_O$ with interpolation parameter $a/k\in(0,1)$. Together
with \eqref{eq:zeroth}, this gives
\begin{align*}
\|A_{\xi}^{a/2}F_\xi\|_2^2
&\le
\bigl(\|A_{\xi}^{k/2}F_\xi\|_2^2\bigr)^{a/k}
\bigl(\|F_\xi\|_2^2\bigr)^{1-a/k}\\
&\le C_a \mathopen{\langle}{\xi}\mathclose{\rangle}^{2a-\beta}.
\end{align*}
Substituting this bound into the initial Plancherel identity
proves \eqref{eq:operator-estimate}, since
\[
\int_{\R^d}|\wh f(\xi)|^2
\mathopen{\langle}{\xi}\mathclose{\rangle}^{2a-\beta}\,d\xi
=\|f\|_{H^{a-\beta/2}}^2.
\]
Taking $a=\beta/2$ and extending from Schwartz functions by density gives
the stated bounded map from $L^2$ to $H^{\beta/2}$. This completes
the proof.
\end{proof}

\subsection{Pairing with the parameter measure}
The full Sobolev estimate can now be paired with a measure on the joint parameter set. We first prove the negative-order Sobolev bound supplied by its Frostman condition.

\begin{proof}[Proof of Lemma~\ref{lem:negative}]
Choose a finite smooth partition of unity near $\spt\tau$,
with nonnegative cutoffs supported in coordinate systems. The coordinate
maps and their inverses are Lipschitz on these supports, so each
localized measure remains $t$-Frostman. The definition of Sobolev
regularity in systems therefore reduces the claim to a compactly
supported measure on $\R^N$.

We now estimate its $H^{-a}$ norm using the Bessel kernel of
$(1-\Delta)^{-a}$. There are three cases, according to the singularity
of this kernel at the diagonal. If $0<2a<N$, the kernel is bounded by
$C|x-y|^{2a-N}$ near the diagonal and decays at infinity. Hence,
\[
\|\tau\|_{H^{-a}}^2
\le C\left(\tau(\R^N)^2+
\iint_{|{x}-{y}|<1}|{x}-{y}|^{-(N-2a)}\,d\tau({x})d\tau({y})\right).
\]
Splitting $0<|x-y|<1$ into dyadic annuli and applying the
$t$-Frostman bound gives, uniformly in $x$,
\begin{align*}
\int_{|x-y|<1}|x-y|^{-(N-2a)}\,d\tau(y)
&\le C\sum_{j\ge0}2^{j(N-2a)}\tau(B(x,2^{-j}))\\
&\le C\sum_{j\ge0}2^{-j(t-N+2a)}<\infty.
\end{align*}
Here $t-N+2a>0$ by hypothesis; in particular, $t>0$, so
$\tau$ has no atoms and the diagonal contributes no mass.

If $2a=N$, the kernel is bounded near the diagonal by
$C\log(2/|x-y|)$. The hypothesis now gives $t>0$, and the same
annular decomposition yields
\[
\int_{|x-y|<1}\log\frac{2}{|x-y|}\,d\tau(y)
\lesssim\sum_{j\ge0}(j+1)2^{-jt}<\infty.
\]
If $2a>N$, the kernel is bounded near the diagonal, so finiteness of
$\tau$ suffices. These cases cover every $a>0$ in the statement.
Each localized measure belongs to $H^{-a}(\R^N)$, and the finite
partition of unity gives $\tau\in H^{-a}(M)$. This completes the proof.
\end{proof}

The next proposition transfers the averaging estimate to the image measure. It uses the joint regularity of $\tau$ and requires no estimates for its individual translation fibres.
\begin{proposition}\label{prop:measure}
Let $\mu$ be a compactly supported finite positive measure satisfying \eqref{eq:spherical}, and let $\tau$ be a compactly supported finite positive measure on $ \Isom(\R^d)$. Define
\[
\nu=\bigl((x,g,z)\mapsto gx+z\bigr)_\#(\mu\times\tau).
\]
If $\tau\in H^{-a}( \Isom(\R^d))$ for some $a\ge0$, then
\begin{equation}\label{eq:image-sobolev}
\nu\in H^{\beta/2-a}(\R^d),\qquad
\|\nu\|_{H^{\beta/2-a}}
\le C_a\|\tau\|_{H^{-a}}.
\end{equation}
In particular, $\tau\in H^{-\beta/2}$ implies that $\nu$ has an $L^2$ density.
\end{proposition}

\begin{proof}
Let $f$ be a Schwartz function. The function $T_\mu f$ is
smooth and rapidly decreasing in $z$, uniformly for $g\in O(d)$,
because $\mu$ has compact support. Thus, its Sobolev pairing with
$\tau$ agrees with integration. Fubini's theorem, Sobolev duality,
and Proposition~\ref{prop:operator} at order $a$ give
\begin{align*}
\left|\int f\,d\nu\right|
&=|\langle\tau,T_\mu f\rangle|\\
&\le \|\tau\|_{H^{-a}}\|T_\mu f\|_{H^a}\\
&\le C_a\|\tau\|_{H^{-a}}\|f\|_{H^{a-\beta/2}}.
\end{align*}
Hence, $\nu$ defines a bounded functional on
$H^{a-\beta/2}(\R^d)$. Identifying its dual with
$H^{\beta/2-a}(\R^d)$ proves \eqref{eq:image-sobolev}.
For $a=\beta/2$, this space is $H^0=L^2$, which proves the final
assertion. This completes the proof.
\end{proof}

\begin{proof}[Proof of Theorem~\ref{thm:main}]
For part~(i), it is enough to construct a probability measure
supported on a compact subset of $\Theta(E)$ whose weighted Fourier
integral is finite at each exponent strictly below the claimed bound.
Fix
\[
{u}_{0}<\dimSF E+\dimH\Theta-m.
\]
\begingroup
Choose exponents $\beta,t\ge0$ with $u_0<\beta+t-m$.
When the corresponding dimension is positive, choose
$\beta<\dimSF E$ and $t<\dimH\Theta$. The definition of spherical
Fourier dimension and Frostman's lemma give probability measures
$\mu\in\mathcal M(E)$ and $\tau\in\mathcal M(\Theta)$ satisfying
\eqref{eq:spherical} and the $t$-Frostman bound, respectively.
If either dimension is zero, use exponent zero and a point mass on
that set; these satisfy the required bounds.
\endgroup
Since $n=m+d$, the inequality $u_0<\beta+t-m$ ensures that
the following interval is nonempty. Choose $a>0$ such that
\[
\frac{n-t}{2}<a<\frac{\beta+d-{u}_{0}}{2}.
\]
The lower bound on $a$ allows us to apply
Lemma~\ref{lem:negative} on the $n$-dimensional motion space.
Thus, $\tau\in H^{-a}$, and Proposition~\ref{prop:measure} places the
image probability measure $\nu$ in $H^{\beta/2-a}(\R^d)$.
The upper bound on $a$ gives $\beta/2-a>(u_0-d)/2$, so
\[
\int|\wh\nu(\xi)|^2(1+|\xi|)^{{u}_{0}-d}\,d\xi<\infty.
\]
Its support is a compact subset of $\Theta(E)$. Taking the supremum over ${u}_{0}$ in \eqref{eq:sobolev-dimension} proves part~(i).

For part~(ii), choose $\beta$ and $t$ as above with
$\beta+t>n$. Since $t\le n$, we have $\beta>0$.
Lemma~\ref{lem:negative} therefore applies with $a=\beta/2$ and gives $\tau\in H^{-\beta/2}(\Isom(\R^d))$. Proposition~\ref{prop:measure} yields an $L^2$ density for $\nu$. Since $\nu$ has mass one, the Cauchy--Schwarz inequality gives
\[
1\le\Leb^d(\spt\nu)\|\nu\|_2^2.
\]
Thus, $\Leb^d(\Theta(E))>0$.

For part~(iii), choose $u_0$ with
$u<u_0<\dimSF E+\dimH\Theta-m$. The preceding construction gives
$\nu\in\mathcal M(\Theta(E))$ with finite weighted Fourier integral
at exponent $u_0$, and hence at exponent $u$. Since $\wh\nu$ is
bounded near the origin and $0<u<d$, the energy identity
\eqref{eq:riesz-energy} gives $I_u(\nu)<\infty$.
Therefore, $\dimH\Theta(E)\ge u$. This completes the proof.
\end{proof}

\begin{proof}[Proof of Corollary~\ref{cor:hausdorff}]
By \eqref{eq:sf-lower},
{
\[
\dimSF E+\dimH\Theta
\ge\frac{d-1}{d}\dimH E+\dimH\Theta.
\]
}
Thus, the hypothesis of part~(i) implies the hypothesis of
Theorem~\ref{thm:main}(ii), and the hypothesis of part~(ii) implies
that of Theorem~\ref{thm:main}(iii), with the same $u$.
Applying these two parts of the theorem proves the corollary.
\end{proof}

\begin{remark}\label{rem:endpoint-energy}
Theorem~\ref{thm:main} also gives positive Lebesgue measure at the
critical dimension under an additional energy assumption.
Indeed, if $\mu$ satisfies (1.2) with $0<\beta\le d$
and $\Theta$ supports a probability measure with finite
$(n-\beta)$-energy, then Proposition~\ref{prop:measure}
implies that $\Theta(E)$ has positive Lebesgue measure.
However, the condition $\dimH\Theta=n-\beta$ alone
does not guarantee finite $(n-\beta)$-energy.
\end{remark}

\section{Proof of Theorem~\ref{thm:similarity}}\label{sec:similarity-proof}
In this section, we prove Theorem~\ref{thm:similarity}. Averaging over rotations and scales replaces spherical Fourier decay by the Fourier estimate on balls for a Frostman source measure. We first prove that estimate and then use Sobolev duality on the similarity parameters.

\begin{proof}[Proof of Lemma~\ref{lem:ball}]
\begingroup
We first bound the indicator function of the frequency ball by a Gaussian.
Fourier inversion then gives, for $R\ge1$,
\begin{align*}
\int_{|\xi|\le R}|\wh{w\mu}(\xi)|^2\,d\xi
&\le e^\pi\int e^{-\pi|\xi|^2/R^2}|\wh{w\mu}(\xi)|^2\,d\xi\\
&=e^\pi R^d\iint e^{-\pi R^2|x-y|^2}w(x)\overline{w(y)}\,d\mu(x)d\mu(y)\\
&\le e^\pi R^d\|w\|_\infty^2
\iint e^{-\pi R^2|x-y|^2}\,d\mu(x)d\mu(y).
\end{align*}
Now split the inner integral
into $B(x,R^{-1})$ and the annuli
$2^jR^{-1}<|x-y|\le2^{j+1}R^{-1}$, $j\ge0$. Uniformly in $x$, this gives
\begin{align*}
\int e^{-\pi R^2|x-y|^2}\,d\mu(y)
&\le \mu(B(x,R^{-1}))+
\sum_{j\ge0}e^{-\pi4^j}\mu(B(x,2^{j+1}R^{-1}))\\
&\lesssim R^{-s}\left(1+\sum_{j\ge0}2^{(j+1)s}e^{-\pi4^j}\right)
\lesssim R^{-s}.
\end{align*}
Integrating this uniform estimate in $x$ and substituting
into the Gaussian bound proves the asserted estimate. The argument
also covers $s=0$ and $s=d$. This completes the proof.
\endgroup
\end{proof}

We now apply this estimate to the similarity averaging operator.

\begin{proof}[Proof of Theorem~\ref{thm:similarity}]
Choose $s<\dimH E$ and $t<\dimH\Omega$, and let
$\mu\in\mathcal M(E)$ and $\tau\in\mathcal M(\Omega)$
be $s$- and $t$-Frostman probability measures, respectively.
We first prove the estimates for these exponents and then
let $s$ and $t$ approach the corresponding dimensions.

Since $\spt\tau$ is compact, its scale projection lies in the
interior of an interval $[r_0,r_1]\subset(0,\infty)$. Define
\[
\mathcal T_\mu f(g,z,r)=\int f(rgx+z)\,d\mu(x),
\]
and choose $\chi\in C_c^\infty((r_0,r_1))$ equal to one on that scale
projection. The estimate we need is
\begin{equation}\label{eq:similarity-operator}
\|\chi\mathcal T_\mu f\|_{H^a(\Sim(\R^d))}
\lesssim_a\|f\|_{H^{a-s/2}(\R^d)},\qquad a\ge0.
\end{equation}
We first explain its consequences, assuming it for the moment.
Let $\nu$ be the image of $\mu\times\tau$ under
$(x,g,z,r)\mapsto rgx+z$. This is a probability measure supported on a
compact subset of $\Omega(E)$. The parameter space has dimension $n+1$,
so Lemma~\ref{lem:negative} gives $\tau\in H^{-a}$ whenever
$a>0$ and $2a>n+1-t$. Since $\chi=1$ on the scale projection of
$\spt\tau$, Fubini's theorem and \eqref{eq:similarity-operator} give
\[
\left|\int f\,d\nu\right|
=|\langle\tau,\chi\mathcal T_\mu f\rangle|
\lesssim_a\|\tau\|_{H^{-a}}\|f\|_{H^{a-s/2}}.
\]
Sobolev duality therefore yields $\nu\in H^{s/2-a}(\R^d)$.

For part~(i), fix $u_0<s+t-(m+1)$. Since $n=m+d$, this
inequality allows us to choose
\[
\frac{n+1-t}{2}<a<\frac{s+d-u_0}{2}.
\]
Then $s/2-a>(u_0-d)/2$, so the weighted Fourier integral of $\nu$ at
exponent $u_0$ is finite. Taking the supremum over these $u_0$ and
letting $s,t$ approach the dimensions proves part~(i), subject to
\eqref{eq:similarity-operator}. For part~(ii), choose $s,t$ with
$s+t>n+1$. This forces $s>0$, and $a=s/2$ gives an $L^2$ image density.
Since $\nu$ has mass one, its support has positive Lebesgue measure.
For part~(iii), choose $s,t$ with $u<s+t-(m+1)$ and then choose $a$
so that $s/2-a>(u-d)/2$. The energy identity \eqref{eq:riesz-energy},
using $0<u<d$ and boundedness of $\wh\nu$ near the origin, gives
$I_u(\nu)<\infty$ and hence $\dimH\Omega(E)\ge u$.

It remains to prove \eqref{eq:similarity-operator} for the
fixed $s$-Frostman measure $\mu$. To control derivatives in the rotation
and scale variables, we first estimate the rotation-scale average of
$\wh{w\mu}$ for a fixed smooth compactly supported weight $w$.
For $|\xi|\ge1$, polar coordinates and Lemma~\ref{lem:ball} give
\begin{align}
\int_{r_0}^{r_1}\int_{O(d)}|\wh{w\mu}(rg^{-1}\xi)|^2\,dg\,dr
&=|\xi|^{-1}\int_{r_0|\xi|}^{r_1|\xi|}\int_{S^{d-1}}
|\wh{w\mu}(\rho\omega)|^2\,d\sigma_{d-1}(\omega)\,d\rho\notag\\
&\lesssim |\xi|^{-d}\int_{r_0|\xi|\le|{\zeta}|\le r_1|\xi|}|\wh{w\mu}({\zeta})|^2\,d{\zeta}\notag\\
&\lesssim |\xi|^{-s}.\label{eq:annular-average}
\end{align}
For the last step, apply Lemma~\ref{lem:ball} with radius
$\max\{1,r_1\}|\xi|$. The constants may depend on $r_0,r_1$ and $w$.
For $|\xi|<1$, the same average is bounded by finiteness of $\mu$.

For integer order $a=k\ge0$, take the Fourier transform in $z$:
\[
\mathcal F_z(\chi\mathcal T_\mu f)(g,\xi,r)
=\chi(r)\wh f(\xi)\wh\mu(-rg^{-1}\xi).
\]
Each derivative in $g$ or $r$ introduces a polynomial weight in $x$
and at most one frequency factor; derivatives falling on $\chi$ have
bounded coefficients. Each derivative in $z$ also introduces one
frequency factor. Thus, derivatives of total order at most $k$ give
finite sums of terms of the form
\[
\xi^\alpha\wh f(\xi)\wh{w\mu}(-rg^{-1}\xi),\qquad |\alpha|\le k,
\]
with bounded coefficients supported in $r\in[r_0,r_1]$.
Here each polynomial weight on $\spt\mu$ is extended to a fixed smooth
compactly supported function $w$. Use finitely many local rotation coordinates and smooth cutoffs
to express the integer Sobolev norm in terms of these derivatives. Then \eqref{eq:annular-average}, the bounded-frequency
estimate, and Plancherel give
\[
\|\chi\mathcal T_\mu f\|_{H^k}^2
\lesssim_k\int|\wh f(\xi)|^2(1+4\pi^2|\xi|^2)^{k-s/2}\,d\xi.
\]
Interpolation between the bounded maps at consecutive
integer orders gives \eqref{eq:similarity-operator} for every real
$a\ge0$. The conditional deductions at the start of the proof now
establish all three assertions. This completes the proof.
\end{proof}

\Needspace{12\baselineskip}
\section{Proof of Theorem~\ref{app-thm-stationary}}\label{sec:new-applications}

In this section, we prove  Theorem~\ref{app-thm-stationary} by extending the similarity averaging estimate to distributions of arbitrary Sobolev order and using stationarity to obtain smooth densities. We then construct a stationary probability measure supported on $E$ to prove that $E$ has nonempty interior, as asserted in Corollary~\ref{app-cor-self-embeddings}.  Throughout this section, we allow $d\ge1$.

For a compactly supported distribution $v$ and a compactly supported finite signed measure $\eta$ on $\Sim(\R^d)$, define $P_\eta v$ by
\[
 \langle P_\eta {v},f\rangle
 =\int\langle {v},f\circ\theta\rangle\,d\eta(\theta),
 \qquad f\in C_c^\infty(\R^d).
\]

\Needspace{13\baselineskip}
\begin{lemma}\label{app-lem-all-orders}
Fix compact sets $K_{x}\subset\R^d$ and $K_{\mathrm{par}}\subset\Sim(\R^d)$.  For every
${\ell}\in\R$ and $a\ge0$,
\begin{equation}\label{app-eq-all-orders}
 \|P_\eta {v}\|_{H^{{\ell}+d/2-a}(\R^d)}
 \lesssim_{K_{x},K_{\mathrm{par}},{\ell},a}
 \|\eta\|_{H^{-a}(\Sim(\R^d))}\|{v}\|_{H^{\ell}(\R^d)}
\end{equation}
whenever $v\in H^\ell(\R^d)$ is supported on $K_x$ and $\eta\in H^{-a}(\Sim(\R^d))$ is a finite signed measure supported on $K_{\mathrm{par}}$. In particular, if $\tau$ is a $t$-Frostman probability measure supported on $K_{\mathrm{par}}$, with $t>m+1$, then
\begin{equation}\label{app-eq-sobolev-gain}
 \|P_\tau {v}\|_{H^{{\ell}+{\varepsilon}_{0}}(\R^d)}
 \lesssim\|{v}\|_{H^{\ell}(\R^d)},
 \qquad 0<{\varepsilon}_{0}<\frac{t-m-1}{2}.
\end{equation}
\end{lemma}

\begin{proof}
Choose $0<r_0<r_1$ so that the scale projection of
$K_{\mathrm{par}}$ lies in $(r_0,r_1)$, and let
$\chi\in C_c^\infty((r_0,r_1))$ equal one near that projection.
For a Schwartz function $f$ on $\R^d$, set
\[
 \mathcal T_v f(g,z,r)=\langle v,x\mapsto f(rgx+z)\rangle.
\]
We will prove
\begin{equation}\label{app-eq-distribution-operator}
 \|\chi\mathcal T_v f\|_{H^a(\Sim(\R^d))}
 \lesssim\|v\|_{H^\ell}\|f\|_{H^{a-d/2-\ell}},
 \qquad a\ge0.
\end{equation}
This estimate suffices for \eqref{app-eq-all-orders}. Indeed, since
$\chi=1$ near the scale projection of $\spt\eta$, it would give
\[
 |\langle P_\eta v,f\rangle|
 =|\langle\eta,\chi\mathcal T_v f\rangle|
 \lesssim\|\eta\|_{H^{-a}}\|v\|_{H^\ell}
 \|f\|_{H^{a-d/2-\ell}}
\]
by Sobolev duality on the parameter space. Duality on $\R^d$ would then
place $P_\eta v$ in $H^{\ell+d/2-a}$ with the required norm bound.

To prove \eqref{app-eq-distribution-operator}, we first estimate
the rotation and scale averages that arise after taking the Fourier
transform in $z$. For every $w\in C_c^\infty(\R^d)$ and $|\xi|\ge1$,
Haar averaging and polar coordinates give
\begin{align}
 \int_{r_0}^{r_1}\int_{O(d)}
 |\wh{w{v}}(rg^{-1}\xi)|^2\,dg\,dr
 &\lesssim {|}{\xi}{|}^{-d}
   \int_{r_0{|}{\xi}{|}\le|\zeta|\le r_1{|}{\xi}{|}}|\wh{w{v}}(\zeta)|^2\,d\zeta\notag\\
 &\lesssim {|}{\xi}{|}^{-d-2{\ell}}\|w{v}\|_{H^{\ell}}^2
 \lesssim_w {|}{\xi}{|}^{-d-2{\ell}}\|{v}\|_{H^{\ell}}^2.
 \label{app-eq-annular-sobolev}
\end{align}
On this annulus, $(1+4\pi^2|\zeta|^2)^\ell$ is comparable to
$|\xi|^{2\ell}$ for every fixed real $\ell$, which proves the second
inequality. The last inequality uses boundedness of multiplication by
$w$ on $H^\ell$. For $|\xi|\le1$, Sobolev duality bounds
$|\wh{wv}(rg^{-1}\xi)|$ by $C_w\|v\|_{H^\ell}$: the corresponding
smooth test functions have uniformly bounded $H^{-\ell}$ norms.
For $d=1$, the same calculation uses the two elements of
$O(1)=\{-1,1\}$ and the scale interval.

The change of variables in the translation variable gives
\[
 \mathcal {F}_{z}\mathopen{(}\mathcal {T}_{v} {f}\mathclose{)}\mathopen{(}{g}\mathpunct{,}{\xi}\mathpunct{,}{r}\mathclose{)}
 \mathrel{=}\wh f(\xi)\wh {v}(-rg^{-1}\xi){.}
\]
Each rotation or scale derivative gives a polynomial
weight in $x$ and at most one factor of $\xi$. Translation derivatives
also give one factor of $\xi$, while derivatives of $\chi$ give
bounded coefficients supported in $(r_0,r_1)$. Since $v$ is supported
on $K_x$, the polynomial weights may be multiplied by a fixed smooth
cutoff equal to one near $K_x$. We can therefore apply
\eqref{app-eq-annular-sobolev} to every resulting weight.
Using finitely many rotation coordinates, Plancherel's theorem in $z$, and
the bounded-frequency estimate, we obtain, for every integer $k\ge0$,
\[
 \|\chi\mathcal T_v f\|_{H^k}^2
 \lesssim \|v\|_{H^\ell}^2
 \int_{\R^d}|\wh f(\xi)|^2
 (1+|\xi|^2)^{k-d/2-\ell}\,d\xi.
\]
For fixed $v$, interpolation of the linear map
$f\mapsto\chi\mathcal T_v f$ between consecutive integer orders proves
\eqref{app-eq-distribution-operator} for all $a\ge0$. The initial
Sobolev duality argument now proves \eqref{app-eq-all-orders}.

For the final assertion, the parameter manifold has dimension
$n+1=m+d+1$. A nonzero $t$-Frostman measure on this manifold has
$t\le n+1$. Thus, if
$0<\varepsilon_0<(t-m-1)/2$, then
$a=d/2-\varepsilon_0>0$ and
\[
 2a=d-2\varepsilon_0>n+1-t.
\]
Lemma~\ref{lem:negative}, applied to $\tau$ on $\Sim(\R^d)$, gives
$\tau\in H^{-a}$. Substituting this value of $a$ into
\eqref{app-eq-all-orders} proves \eqref{app-eq-sobolev-gain} and
completes the proof.
\end{proof}

\begin{proof}[Proof of Theorem~\ref{app-thm-stationary}]
We first use stationarity to iterate the positive Sobolev gain
under the Frostman hypothesis. We then treat the critical energy
condition by a smooth approximation of the parameter law.

First suppose that $\tau$ is $t$-Frostman with $t>m+1$.  A compactly
supported probability measure $\mu$ belongs to $H^{{\ell}_0}$ for every
${\ell}_0<-d/2$.  Fix $0<{\varepsilon}_{0}<(t-m-1)/2$.  Stationarity and
\eqref{app-eq-sobolev-gain} imply successively that
\[
 \mu\in H^{{\ell}_0+j{\varepsilon}_{0}}(\R^d),\qquad j=0,1,2,\ldots.
\]
At every step, the distribution is $\mu$, so its support
remains in the same compact set. For any integer $k\ge0$, finitely many
iterations give $\mu\in H^\ell$ for some $\ell>k+d/2$. Sobolev
embedding then gives a $C^k$ density. Since $k$ is arbitrary, the
density is $C^\infty$.

Now assume only \eqref{app-eq-critical-law}.  The case $a=d/2$ of
\eqref{app-eq-all-orders} gives
\begin{equation}\label{app-eq-critical-operator}
 \|P_\eta {v}\|_{H^{\ell}}
 \lesssim_{\ell}\|\eta\|_{H^{-d/2}}\|{v}\|_{H^{\ell}},
 \qquad {\ell}\in\R,
\end{equation}
with constants uniform on fixed compact source and parameter sets.
We will use this estimate to make the nonsmooth part of the
parameter law define a contraction at two prescribed Sobolev orders.
Choose $\psi\in C_c^\infty(\R^d)$ equal to one near $\spt\mu$.
By a finite partition of unity in parameter coordinates and convolution
with smooth approximate identities, choose smooth compactly supported
densities $\tau_\varepsilon$ such that
\[
 \|\tau-\tau_\varepsilon\|_{H^{-d/2}(\Sim(\R^d))}\longrightarrow0.
\]
Their supports lie in one enlarged compact parameter set whose scales
are bounded away from zero. Set
\[
 \eta_\varepsilon=\tau-\tau_\varepsilon,
 \qquad Q_\varepsilon v=\psi P_{\eta_\varepsilon}(\psi v),
 \qquad S_\varepsilon v=\psi P_{\tau_\varepsilon}(\psi v).
\]
The input cutoff places $\psi v$ in the fixed compact source set
$\spt\psi$, so \eqref{app-eq-critical-operator} applies uniformly in
$\varepsilon$. Boundedness of multiplication by $\psi$ on $H^\ell$
then gives
\[
 \|Q_\varepsilon\|_{H^\ell\to H^\ell}
 \lesssim_\ell\|\eta_\varepsilon\|_{H^{-d/2}}
 \longrightarrow0
\]
for every fixed $\ell\in\R$.
Writing $d\tau_\varepsilon=h_\varepsilon(g,z,r)\,dg\,dz\,dr$,
the operator $S_\varepsilon$ has the smooth compactly supported kernel
\[
 K_\varepsilon(y,x)=\psi(y)\psi(x)
   \int h_\varepsilon(g,y-rgx,r)\,dg\,dr.
\]
In particular, $S_\varepsilon\mu$ is smooth and compactly supported.

Fix $\ell_0<-d/2$ and a target order $\ell_1>\ell_0$.
Choose $\varepsilon$ small enough that $Q_\varepsilon$ has norm less
than $1/2$ on both $H^{\ell_0}$ and $H^{\ell_1}$.
Since $\psi\mu=\mu$ and $P_\tau\mu=\mu$, we have
\[
 \mu=Q_\varepsilon\mu+S_\varepsilon\mu,
 \qquad
 (\operatorname{Id}-Q_\varepsilon)\mu=S_\varepsilon\mu.
\]
The Neumann series
\[
 \sum_{j=0}^{\infty}Q_\varepsilon^jS_\varepsilon\mu
\]
converges in $H^{\ell_1}$, because
\[
 \|Q_\varepsilon^jS_\varepsilon\mu\|_{H^{\ell_1}}
 \le2^{-j}\|S_\varepsilon\mu\|_{H^{\ell_1}}.
\]
Its sum also belongs to $H^{\ell_0}$, by the continuous inclusion
$H^{\ell_1}\subset H^{\ell_0}$, and solves the same equation as $\mu$.
The operator $\operatorname{Id}-Q_\varepsilon$ is invertible on
$H^{\ell_0}$, so uniqueness there identifies the sum with $\mu$.
Hence, $\mu\in H^{\ell_1}$. Since the target order $\ell_1$ was
arbitrary, Sobolev embedding gives a $C^\infty$ density. This completes
the proof.
\end{proof}

\begin{proof}[Proof of Corollary~\ref{app-cor-self-embeddings}]
We construct a stationary probability measure supported on ${E}$.
Theorem~\ref{app-thm-stationary} gives a smooth density, which must be
positive on some ball.

We first choose a parameter to which that theorem applies.
On the parameter manifold, of dimension $n+1$, the Bessel kernel for
the squared $H^{-d/2}$ norm has local singularity of order
$(n+1)-d=m+1$. For a positive compactly supported measure, the kernel
away from the diagonal contributes a finite quantity. Thus, by the
local-coordinate argument used in Lemma~\ref{lem:negative}, finite
$(m+1)$-energy is equivalent to \eqref{app-eq-critical-law}.

Choose a probability measure $\tau$ supported on $\Phi$ with
finite $(m+1)$-energy and choose
$m+1<t<\dimH\Phi$ and a $t$-Frostman probability measure on $\Phi$. Since every $\theta\in\Phi$ maps $E$ into itself, $P_\tau$ maps
probabilities on $E$ into probabilities on $E$.  It is weakly continuous:
for $f\in C(E)$, the function
$x\mapsto\int f(\theta x)\,d\tau(\theta)$ is continuous by uniform
continuity on $\Phi\times E$. Let $\mu_0$ be a probability measure on $E$,
and set
\[
 {\mu}_{J}\mathrel{=}\frac1{J}\sum_{j=0}^{{J}-1}P_\tau^j\mu_0.
\]
These averages are approximately stationary, since the sum
telescopes to give
\[
 P_\tau{\mu}_{J}\mathbin{-}{\mu}_{J}\mathrel{=}\frac{{P}_{\tau}^{J}{\mu}_{0}\mathbin{-}{\mu}_{0}}{{J}}\mathpunct{,}
 \qquad
 {|}\mathopen{\langle} {P}_{\tau}{\mu}_{J}\mathbin{-}{\mu}_{J}\mathpunct{,}{f}\mathclose{\rangle}{|}
 \mathrel{\le}\frac{{2}{\|}{f}{\|}_{\infty}}{{J}}{.}
\]
Compactness of $E$ gives a weakly convergent subsequence of
$(\mu_J)$ whose limit $\mu$ is a probability measure on $E$.
The above estimate and weak continuity of $P_\tau$ imply
$P_\tau\mu=\mu$. Theorem~\ref{app-thm-stationary}, applied with the
parameter law $\tau$, gives a smooth nonnegative density of integral
one. This density is positive on some open ball. That ball lies in
$E$: otherwise, since $E$ is closed, its intersection with the
complement of $E$ would contain an open ball of positive $\mu$-mass,
contrary to $\spt\mu\subset E$. Hence, $E$ has nonempty interior.

It remains to deduce the assertions about the full family
$\Sigma_E$. This family is closed in $\Sim(\R^d)$: if
$\theta_j\to\theta$ and $\theta_j(E)\subset E$, then
$\theta_j(x)\to\theta(x)\in E$ for every $x\in E$, because $E$ is
closed. For each integer $j\ge1$, let
\[
\Phi_j
=
\Sigma_E\cap
\bigl(O(d)\times B(0,j)\times[1/j,j]\bigr).
\]
Each $\Phi_j$ is compact, and
\[
\Sigma_E=\bigcup_{j=1}^{\infty}\Phi_j.
\]
If $\dimH\Sigma_E>m+1$, then the countable stability
of Hausdorff dimension implies that
$\dimH\Phi_j>m+1$ for some $j$.
By the preceding argument, $E$ must have nonempty interior.

If $\Sigma_E$ supports a compactly supported probability measure
of finite $(m+1)$-energy, apply the first part with $\Phi$ equal to its
support. This again gives nonempty interior. Taking contrapositives
proves both assertions for a set $E$ with empty interior and completes
the proof.
\end{proof}

\begingroup
To see sharpness when $d\ge2$, let $H$ be a linear hyperplane, and let
$E$ be its closed unit ball. Define
\[
 \Phi=\{(g,z,r)\in\Sim(\R^d):gH=H,\ r\in[1/4,1/2],\ z\in H,\ |z|\le1-r\}.
\]
The family $\Phi$ is compact. If $x\in E$ and
$(g,z,r)\in\Phi$, then $rgx+z\in H$ and
$|rgx+z|\le r+|z|\le1$. Thus, every similarity in $\Phi$ maps $E$
into itself. The rotations preserving $H$ form a group of dimension
$(d-1)(d-2)/2$. For each $r\in[1/4,1/2]$, the allowed translations
form a ball with nonempty relative interior in $H$. Consequently,
\[
 \dimH\Phi=\frac{(d-1)(d-2)}2+1+(d-1)=m+1.
\]
For every $y\in E$ and $r\in[1/4,1/2]$,
\[
 y=ry+(1-r)y,\qquad
 (\operatorname{Id},(1-r)y,r)\in\Phi.
\]
This proves the reverse inclusion $E\subset\Phi(E)$, and hence
$\Phi(E)=E$. The map sending a nonempty compact set $A$ to $\Phi(A)$
is a contraction with constant at most $1/2$ in the Hausdorff metric.
Thus, its fixed set $E$ is the attractor of $\Phi$, although $E$ has
empty interior in $\R^d$.
For $d=1$, take $E=\{0\}$ and the family $x\mapsto rx$,
$r\in[1/4,1/2]$. This family has dimension $1=m+1$ and attractor $E$.
Therefore, the dimension threshold cannot be lowered without further
assumptions.
\endgroup

\section{Sharpness examples}\label{sec:examples}
In this section, we show that the strict dimension inequalities in the positive-measure conclusions of Theorems~\ref{thm:main} and~\ref{thm:similarity} and Proposition~\ref{prop:line-area} cannot be relaxed to equality. We also prove that the coefficient $(d-1)/d$ in Corollary~\ref{cor:hausdorff} cannot be increased. Finally, we construct sharp examples for every fixed nonempty finite set and every possible dimension of the motion family.

\begin{proposition}[]\label{prop:sharp}
For every $d\ge2$, there are compact $E\subset\R^d$ and compact $\Theta\subset \Isom(\R^d)$ such that
\[
\dimSF E+\dimH\Theta=n,
\qquad \Leb^d(\Theta(E))=0.
\]
One may take $E$ to be the unit sphere and the rotation set to be all of $ O(d)$.
\end{proposition}

\begin{proof}
Choose a compact set $C\subset\R$ with
$\dimH C=1$ and $\Leb^1(C)=0$. Translating the unit sphere along $Ce_1$
will give a measure zero union whose one-dimensional sections are finite unions
of translates of $C$. Set
\[
E=S^{d-1},\qquad \Theta=O(d)\times(Ce_1).
\]
The normalized surface measure $\mu$ on $S^{d-1}$ satisfies
\[
|\wh\mu(\xi)|\lesssim (1+|\xi|)^{-(d-1)/2}.
\]
Squaring and averaging over directions gives $\dimSF E\ge d-1$.
Together with $\dimSF E\le\dimH E=d-1$, this yields
$\dimSF E=d-1$. Since $O(d)$ is an $m$-dimensional compact manifold,
the product formula with a smooth manifold gives
$\dimH\Theta=m+\dimH C=m+1$. Thus,
$\dimSF E+\dimH\Theta=n$. It remains to verify that the union is of measure zero.
Since rotations preserve $S^{d-1}$, we have
\[
\Theta(E)=S^{d-1}+Ce_1.
\]
Write a point in $\R^d$ as $(x,y)\in\R\times\R^{d-1}$.
For each fixed $y$ with $|y|\le 1$, the section of
$\Theta(E)$ is
\[
\{x:(x,y)\in\Theta(E)\}
=
\bigl(C+\sqrt{1-|y|^2}\bigr)
\cup
\bigl(C-\sqrt{1-|y|^2}\bigr).
\]
Since $\Leb^1(C)=0$, this section has Lebesgue measure
zero. For $|y|>1$, the section is empty.
The set $\Theta(E)$ is compact, so Fubini's theorem applies to its
indicator function and gives
\[
\Leb^d(\Theta(E))=0.
\]
This completes the proof.
\end{proof}

\begin{proposition}[]\label{prop:coefficient}
The coefficient $1/2$ in \eqref{eq:planar-target} cannot be replaced by any larger constant in a theorem for all Borel sets. More generally, the sufficient condition
\[
\frac{d-1}{d}\dimH E+\dimH\Theta>n,
\]
which follows from $\dimSF E\ge\frac{d-1}{d}\dimH E$, fails if $(d-1)/d$ is replaced by any larger coefficient.
\end{proposition}
\begin{proof}
\begingroup
We use the construction in \cite[Example~6.4]{IMPPSS}.
It is enough to construct a compact set $A\subset[0,1]$ with
$\dimH A=1$ and $\Leb^1(A+A)=0$. Indeed, taking its product with a ball
in $\R^{d-1}$ will give a full-dimensional set $E$ whose union under
the motions defined below is contained in
$\R^{d-1}\times(A+A)$. We first construct $A$ and then verify the
dimensions of $E$ and the motion family.

Let $N_k=2^{k^2}$ for $k\ge1$, and set
\[
\mathcal J=\bigcup_{k\ge1}\{N_k+1,\ldots,N_k+k\},
\qquad
A=\left\{\sum_{j=1}^{\infty}\varepsilon_j2^{-j}:
\varepsilon_j\in\{0,1\},\ \varepsilon_j=0\text{ for }j\in\mathcal J\right\}.
\]
The set $A$ is compact. If $N_k\le j<N_{k+1}$, then
$|\mathcal J\cap\{1,\ldots,j\}|\le k(k+1)/2$, whereas
$j\ge2^{k^2}$. Consequently,
$|\mathcal J\cap\{1,\ldots,j\}|=o(j)$ as $j\to\infty$.
Let $\lambda$ be the probability measure obtained by choosing each
unrestricted digit independently with equal probabilities.
Every permitted binary prefix of length $j$ has mass
\[
2^{-j+|\mathcal J\cap\{1,\ldots,j\}|}
\le C_{\varepsilon}2^{-j(1-\varepsilon)},
\qquad \varepsilon>0.
\]
For an interval $I$ with $0<|I|\le1$, choose $j\ge0$ so that
$2^{-j-1}<|I|\le2^{-j}$. Then $I$ meets at most three dyadic intervals
of length $2^{-j}$. Dyadic endpoints have $\lambda$-measure zero,
because infinitely many digits remain free. Applying the prefix
bound to these intervals gives
\[
\lambda(I)\lesssim_{\varepsilon}|I|^{1-\varepsilon}
\qquad(0<|I|\le1),
\]
for every $0<\varepsilon<1$. Frostman's lemma gives
$\dimH A\ge1-\varepsilon$ for each such $\varepsilon$.
Since $A\subset\R$, it follows that $\dimH A=1$.

To prove that $A+A$ is of measure zero, fix $k\ge1$. Truncating the binary
expansion of $x\in A$ after its first $N_k$ digits gives an integer
$j_0$. The next $k$ digits vanish, so the remaining tail satisfies
\[
0\le x-j_0 2^{-N_k}\le 2^{-N_k-k},
\qquad 0\le j_0<2^{N_k}.
\]
For a sum of two points of $A$, the two integers $j_0$ have at
most $2^{N_k+1}$ possible sums, and the combined tail is at most
$2^{1-N_k-k}$. Thus, $A+A$ is covered by at most $2^{N_k+1}$ intervals
of length $2^{1-N_k-k}$. Their total length is at most $2^{2-k}$.
Letting $k\to\infty$ gives $\Leb^1(A+A)=0$.

Set
\[
E=\{(x',x_d)\in\R^{d-1}\times\R:|x'|\le1,\ x_d\in A\},
\]
and let
\[
\Theta=\{(\operatorname{diag}(h,1),z):h\in O(d-1),\ z\in E\}.
\]
The rotations $\operatorname{diag}(h,1)$ preserve $E$, and the
sum of two unit balls in $\R^{d-1}$ is the ball of radius two. Hence,
\[
\Theta(E)=\{(y',y_d)\in\R^{d-1}\times\R:|y'|\le2,\ y_d\in A+A\}
\]
has Lebesgue measure zero by Fubini's theorem. Products with
Euclidean balls and smooth manifolds give
\[
\dimH E=d,
\qquad
\dimH\Theta=\frac{(d-1)(d-2)}2+d=n-d+1.
\]
For every $c>(d-1)/d$, we have
\[
c\dimH E+\dimH\Theta=cd+n-d+1>n.
\]
The same compact sets therefore violate the proposed
positive-measure conclusion for every $c>(d-1)/d$. This proves the
claim in every dimension, including $d=2$.
\endgroup
\end{proof}

\begin{example}[]\label{ex:similarity}
Let $C\subset[0,1]$ be compact, with $\dimH C=1$ and $\Leb^1(C)=0$. Set
\[
E=\{0\},\qquad Z=C\times[0,1]^{d-1},\qquad
\Omega=O(d)\times Z\times[1,2].
\]
Since every similarity sends $0$ to its translation parameter,
$\Omega(E)=Z$. Fubini's theorem gives $\Leb^d(Z)=0$, while the product
dimension formula gives
\[
\dimH E+\dimH\Omega=m+d+1=\frac{d^2+d+2}{2}.
\]
Thus, equality in Theorem~\ref{thm:similarity}(ii) does not suffice. This is the same obstruction as in \cite[Section 6]{IMPPSS}.
\end{example}

\begin{example}[]\label{ex:line-area}
Let $E=[0,1]\times\{0\}$ and let $C\subset[0,1]$ be compact with $\dimH C=1$ and $\Leb^1(C)=0$. For
\[
\Theta=\{(\operatorname{Id},({z}_{1},{z}_{2})):{z}_{1}\in[0,1],\ {z}_{2}\in C\},
\]
one has $\dimH\Theta=1+\dimH C=2$ and
$\Theta(E)=[0,2]\times C$. Fubini's theorem gives
$\Leb^2(\Theta(E))=0$, although $\mathcal H^1(E)=1$.
Thus, equality does not suffice in Proposition~\ref{prop:line-area}.
\end{example}

\begingroup
\begin{proposition}[]\label{prop:finite-source}
Let $d\ge2$, and let $E\subset\R^d$ be any fixed nonempty finite set.
For every $t\in[0,n]$, there is a nonempty compact family
$\Theta\subset\Isom(\R^d)$ such that
\[
\dimH\Theta=t,
\qquad
\dimH\Theta(E)=\max\{0,t-m\}.
\]
In particular, this conclusion holds when $E$ is the vertex set of a
fixed nondegenerate simplex.
\end{proposition}

Since $\dimSF E=0$ for finite $E$, this proposition gives equality in
the Hausdorff lower bound following Theorem~\ref{thm:main} for every
possible dimension of the motion family. Thus, affine spanning of a
finite set does not improve that bound.

\begin{proof}
If $t=0$, we simply take $\Theta$ to be a singleton.
If $t=n$, we take a compact neighborhood of the identity
with nonempty interior. In this case, the image of every
point contains an open ball, because the neighborhood contains
all sufficiently small translations. Thus, its image has dimension
$d=n-m$, as required. Having treated both endpoints, we assume that
$0<t<n$ and set
\[
s=\max\{0,t-m\}.
\]

We will construct $\Theta$ from nested balls.
At each stage, the centers will be chosen to control the image of one
point of $E$. Cycling through these points will give a sequence of
efficient covers for each image. We first describe the coordinates
and reduce the conclusion to the corresponding image upper bounds.
Choose local coordinates $a\mapsto g(a)$ on $O(d)$,
where $a$ ranges over an open subset of $\R^m$..
We identify the motion $(g(a),z)$ with its coordinates $(a,z)$.
For every $x\in E$, define the map
\[
\pi_x(a,z)=g(a)x+z.
\]
We also introduce the change of coordinates
\[
\Psi_x(a,z)=(a,g(a)x+z),
\qquad
\Psi_x^{-1}(a,y)=(a,y-g(a)x).
\]
In these new coordinates, $\pi_x$ is simply the
projection onto the last $d$ coordinates.
Since $E$ is finite, we can choose a fixed compact
coordinate neighborhood on which all the maps
$\Psi_x$ and their inverses have uniform Lipschitz
bounds. Throughout the construction, balls are closed Euclidean balls
in these local coordinates. Since the coordinate map is
bi-Lipschitz on the chosen neighborhood, it preserves
Hausdorff dimension.

For any nonempty compact family $\Theta$ in this coordinate with
$\dimH\Theta=t$, the set $\Psi_x(\Theta)$ is contained in
$[-L,L]^m\times\pi_x(\Theta)$ for some fixed $L$. The product bound
with an $m$-dimensional cube gives
\[
t=\dimH\Psi_x(\Theta)\le m+\dimH\pi_x(\Theta).
\]
Together with nonemptiness, this gives
$\dimH\pi_x(\Theta)\ge\max\{0,t-m\}=s$.
Since $E$ is finite and
$\Theta(E)=\bigcup_{x\in E}\pi_x(\Theta)$, it is therefore enough to
construct a nonempty compact $\Theta$ with $\dimH\Theta=t$ and
$\dimH\pi_x(\Theta)\le s$ for every $x\in E$.

To choose the child centers inside each parent ball, we use a
compact $t$-Ahlfors regular set
$A\subset\R^m\times\R^d$.
If $t\le m$, we take $A$ to be a $t$-Ahlfors
regular subset of $\R^m\times\{0\}$.
If $t>m$, we take
\[
A=[0,1]^m\times C_s,
\]
where $C_s\subset[0,1]^d$ is an $s$-Ahlfors
regular set. Such sets can be constructed using
cubes and self-similar Cantor sets.

\begingroup
Translating a point of $A$ to the origin and rescaling preserves
its Ahlfors regularity and the dimension of its projection onto the
last $d$ coordinates. We make these changes so that $0\in A$ and $A$
lies in a fixed small ball. We then choose a sufficiently small
constant $c_0>0$ so that, for every $x\in E$,
\[
\Psi_x^{-1}\bigl(\Psi_x(\theta)+c_0rA\bigr)
\subset B(\theta,r/4)
\]
whenever $B(\theta,r)$ lies in our chosen coordinate
neighborhood. This follows from the uniform
Lipschitz bounds for $\Psi_x^{-1}$.
Moreover, the set on the left is $t$-Ahlfors regular,
with constants independent of $x$, $\theta$, and $r$
after rescaling by $r$.

The Ahlfors regularity bounds give two estimates that we will
use to count children and bound their local concentration. For
$0<h\ll r$, a maximal $10h$-separated subset of this set contains at
least $c(r/h)^t$ points. Every $10h$-separated subset has at most
$C(\rho/h)^t$ points in a ball of radius $\rho$, for $h\le\rho\le r$.
For the lower bound, maximality gives a covering by balls of radius
$10h$, and each such ball has measure at most a constant times $h^t$.
For the upper bound, use the disjoint balls of radius $5h$ about the
separated points and the lower regularity bound for their measures.

For the image estimate, the same set can be covered after
applying $\pi_x$ by at most $C(r/h)^s$ balls of radius $h$. Indeed,
for $(b,y)\in A$,
\[
\pi_x\!\left(\Psi_x^{-1}\bigl(\Psi_x(\theta)+c_0r(b,y)\bigr)\right)
=\pi_x(\theta)+c_0ry.
\]
The projection of $A$ onto its last $d$ coordinates is $s$-Ahlfors
regular when $s>0$; when $s=0$, the image is a single point.
All constants are independent of the parent ball and of $h$.
\endgroup

List the points of $E$ as $x_1,\ldots,x_{|E|}$ and extend
the list periodically by $x_{j+|E|}=x_j$. We choose scales at which
the logarithmic cost of all previous stages is negligible. This will
allow the image estimate for each point to be proved along its own
subsequence. Choose a sufficiently large integer $J$ and set
\[
r_j=\exp\bigl(-(j+J)!\bigr),\qquad j\ge0.
\]
In particular,
\begin{equation}\label{eq:finite-source-scales}
\frac{\log(1/r_{j-1})}{\log(1/r_j)}\longrightarrow0,
\qquad
\frac{j}{\log(1/r_{j-1})}\longrightarrow0.
\end{equation}
Start with one ball of radius $r_0$ in the interior of the coordinate
neighborhood. \begingroup
In each ball $B(\theta,r_{j-1})$ at stage $j-1$, choose
\[
\left\lfloor c_1(r_{j-1}/r_j)^t\right\rfloor
\]
points in $\Psi_{x_j}^{-1}\bigl(\Psi_{x_j}(\theta)+c_0r_{j-1}A\bigr)$,
separated by
$10r_j$, and take the radius-$r_j$ balls about them.
Choose $c_1>0$ below the constant in the lower bound for the
number of separated points. Then choose $J$ large enough that the
ratio $r_j/r_{j-1}$ is sufficiently small at every stage and every
parent has at least two children. The preceding packing estimate
ensures that the required centers exist. They lie in the central
quarter of the parent, and $r_j<r_{j-1}/4$, so each child ball lies
in its parent.
\endgroup

Let $\mathcal B_j$ be the family of balls at stage $j$. Distinct
balls in $\mathcal B_j$ have centers separated by at least $10r_j$.
Indeed, this holds for children of the same parent by construction;
for different parents it follows from their separation and the
placement of child centers in the central quarter of each parent.
The number of children of each parent is comparable to
$(r_{j-1}/r_j)^t$, with constants independent of $j$. Multiplying
these bounds and telescoping the scale ratios gives, for fixed
constants $0<c<1<C$,
\begin{equation}\label{eq:finite-source-count}
c^j(r_0/r_j)^t
\le |\mathcal B_j|=\prod_{\ell=1}^j
\left\lfloor c_1(r_{\ell-1}/r_\ell)^t\right\rfloor
\le C^j(r_0/r_j)^t.
\end{equation}
\Needspace{6\baselineskip}
The nested intersection
\[
\Theta=\bigcap_{j\ge0}\ \bigcup_{B\in\mathcal B_j}B
\]
is nonempty and compact. We identify it with its image in the motion
group under the fixed coordinate.

We next verify that $\dimH\Theta=t$. Assign mass
$|\mathcal B_j|^{-1}$ to each ball in $\mathcal B_j$.
Every parent at a given stage has the same number of children, so
these assignments are consistent and define a probability measure
$\tau$ on $\Theta$. Suppose that $r_j\le\rho<r_{j-1}$. The separation
of the parent centers implies that a ball of radius $\rho$ meets
only a bounded number of parents at stage $j-1$. Within each parent,
any child ball meeting $B(\theta,\rho)$ has its center in
$B(\theta,\rho+r_j)\subset B(\theta,2\rho)$.
The local counting estimate for child centers,
with a fixed enlargement of the radius, therefore bounds the number
of intersecting children by $C(\rho/r_j)^t$. Hence,
\[
\tau(B(\theta,\rho))
\lesssim\frac{(\rho/r_j)^t}{|\mathcal B_j|}
\lesssim c^{-j}\rho^t.
\]
For every $0<u<t$, \eqref{eq:finite-source-scales} gives
\[
c^{-j}\rho^t
\le c^{-j}r_{j-1}^{t-u}\rho^u
\lesssim_u\rho^u.
\]
For $\rho\ge r_0$, the same bound follows by increasing the constant
and using $\tau(\Theta)=1$. Thus, $\tau$ satisfies every Frostman bound of exponent less than
$t$, and $\dimH\Theta\ge t$. Conversely, for every $v>t$,
the cover by the balls in $\mathcal B_j$ satisfies
\[
|\mathcal B_j|(2r_j)^v\lesssim_v C^j r_j^{v-t}
\longrightarrow0.
\]
Therefore, $\dimH\Theta=t$.

It remains to prove the image upper bounds. Fix $x\in E$ and
consider the infinitely many stages with $x_j=x$.
At such a stage, the image of the child centers in each parent can
be covered by $C(r_{j-1}/r_j)^s$ balls of radius $r_j$.
Every point of $\Theta$ lies within $r_j$ of a child center, and
$\pi_x$ is uniformly Lipschitz. Consequently, $\pi_x(\Theta)$ has
a cover by at most
\[
C|\mathcal B_{j-1}|(r_{j-1}/r_j)^s
\]
balls of radius $Cr_j$. For every $v>s$, the total $v$-dimensional
cost of this cover is at most
\[
C_v C^j r_{j-1}^{s-t}r_j^{v-s}\longrightarrow0,
\]
since \eqref{eq:finite-source-scales} implies
{
\[
\frac{\log\bigl(C^j r_{j-1}^{s-t}r_j^{v-s}\bigr)}{\log(1/r_j)}
=\frac{j\log C+(t-s)\log(1/r_{j-1})}{\log(1/r_j)}-(v-s)
\longrightarrow-(v-s)<0.
\]
}
Thus,
$\dimH\pi_x(\Theta)\le s$ for every $x\in E$.

The lower bound established at the start of the proof now
gives $\dimH\pi_x(\Theta)=s$ for every $x\in E$.
Since Hausdorff dimension is stable under finite unions,
\[
\dimH\Theta(E)
=\dimH\bigcup_{x\in E}\pi_x(\Theta)=s.
\]
This completes the proof.
\end{proof}
\endgroup

\section{Incidence estimates and further applications}\label{sec:applications}
In this section, we prove a quantitative incidence estimate in every dimension and construct examples attaining its three regimes. We also recover Orponen's point-line incidence inequality from Proposition~\ref{prop:operator}. We show that the motions mapping $E$ into a zero set form a family of Hausdorff dimension at most $n-\dimSF E$. For planar patterns contained in a line, we obtain a further dimension bound. If the pattern has positive length, a motion family of dimension greater than two gives a union of positive area.

\subsection{Quantitative incidences under rigid motions}\label{sec:incidences}
\begingroup
We estimate how often a transformed source point lies within distance $\delta$ of a target point.
For compactly supported probability measures $\mu$ on $\R^d$, $\tau$ on $\Isom(\R^d)$, and $\rho$ on $\R^d$, define
\begin{equation}\label{eq:incidence-definition}
\mathcal I_\delta(\mu,\tau,\rho)=
\iiint\mathbf1_{\{|gx+z-y|\le\delta\}}\,d\mu(x)d\tau(g,z)d\rho(y).
\end{equation}

\begin{theorem}\label{thm:incidence}
Let $d\ge2$, and let $\mu$, $\tau$, and $\rho$ be compactly supported probability measures on $\R^d$, $\Isom(\R^d)$, and $\R^d$, respectively. Suppose that $\mu$ satisfies \eqref{eq:spherical} for $0\le\beta\le d$, $\tau$ is $t$-Frostman, and $\rho$ is $b$-Frostman, where $0\le t\le n$ and $0\le b\le d$. For $0<\delta\le1/2$,
\begin{equation}\label{eq:incidence}
\mathcal I_\delta(\mu,\tau,\rho)\lesssim
\begin{cases}
\delta^{(\beta+b+t-m)/2},&\beta+b+t<n+d,\\
\delta^d\log(2/\delta),&\beta+b+t=n+d,\\
\delta^d,&\beta+b+t>n+d.
\end{cases}
\end{equation}
The implicit constant may depend on $d$, the spherical-decay constant, the Frostman constants and exponents, and a bound for $\spt\mu$. In every dimension $d\ge2$, each of the three orders is attained by suitable measures. In particular, the critical logarithmic factor cannot in general be removed.
\end{theorem}

The sharpness examples form a one-parameter family. For every $d\ge2$ and $0<\alpha<1$, Section~\ref{subsec:incidence-sharpness} constructs measures with
\[
(\beta,t,b)=(d-1,m+1+\alpha,d-1+\alpha)
\]
such that, as $\delta\downarrow0$,
\begin{equation}\label{eq:incidence-sharp-family}
\mathcal I_\delta(\mu,\tau,\rho)\sim
\begin{cases}
\delta^{d-1/2+\alpha},&0<\alpha<1/2,\\
\delta^d\log(1/\delta),&\alpha=1/2,\\
\delta^d,&1/2<\alpha<1.
\end{cases}
\end{equation}
The source measure is the product of normalized arc length on a circle and a smooth probability density in the remaining $d-2$ directions. Tangency to a hyperplane produces the same one-dimensional kernel in every dimension. Thus, the bounds are sharp along this family; we do not assert optimality for every admissible triple $(\beta,t,b)$.

The incidence estimate should be compared with the basic bounds
\begin{equation}\label{eq:incidence-basic}
\mathcal I_\delta(\mu,\tau,\rho)
\lesssim\min\bigl\{\delta^b,\delta^{t-m},\delta^{(\beta+b)/2}\bigr\}.
\end{equation}
The target Frostman condition gives the first bound. Covering the $m$-dimensional family of motions sending one fixed point to another gives the second. A Fourier estimate for each fixed motion gives the third; see Lemma~\ref{lem:incidence-basic}. In the regime $\beta+b+t<n+d$, Theorem~\ref{thm:incidence} strictly improves all three bounds precisely when
\[
t>m,\qquad \beta+t>b+m,\qquad \beta+b>t-m.
\]
For other parameter values, one takes the better of \eqref{eq:incidence} and \eqref{eq:incidence-basic}. In particular, if $b=d$ or $t=n$, the basic bounds already give $\delta^d$, so the logarithm is not needed for those critical triples.
\endgroup

We first prove a dyadic Fourier estimate for the image measure, then pair it with the target measure to prove Theorem~\ref{thm:incidence}. We also verify the elementary comparison bounds and construct the sharpness examples. Appendix~\ref{app:planar-fourier} gives an alternative angular Fourier proof in the plane.

\subsubsection{A dyadic estimate in all dimensions}
\begingroup
\begin{proposition}\label{prop:motion-shell}
Let $\mu$ be a compactly supported probability measure on $\R^d$ satisfying \eqref{eq:spherical} for $0\le\beta\le d$, and let $\tau$ be a compactly supported $t$-Frostman probability measure on $\Isom(\R^d)$, where $0\le t\le n$. Let $\nu$ be the image of $\mu\times\tau$ under $(x,g,z)\mapsto gx+z$. Then, for $R\ge1$,
\begin{equation}\label{eq:motion-shell}
\int_{R\le|\xi|\le2R}|\wh\nu(\xi)|^2\,d\xi
\lesssim R^{n-t-\beta}.
\end{equation}
The constant has the same permitted dependencies as in Theorem~\ref{thm:incidence}.
\end{proposition}

\begin{proof}
By $L^2$ duality, it suffices to prove
\[
\left|\int f\,d\nu\right|
\lesssim R^{(n-t-\beta)/2}\|f\|_2
\]
for every Schwartz function $f$ with Fourier support in
$R\le|\xi|\le2R$. Fix such an $f$. We apply weighted Cauchy--Schwarz in local rotation coordinates,
using the Frostman condition for one factor and
Proposition~\ref{prop:operator} for the other.
We only localize the rotation variable, so the constants
do not depend on the translation support of $\tau$.

Choose a finite smooth partition of unity $\{\chi_i\}$ on $O(d)$, with each $\chi_i$ nonnegative and compactly supported in a coordinate system $\kappa_i:U_i\to V_i\subset\R^m$. Choose $\psi_i\in C_c^\infty(U_i)$ equal to one near $\spt\chi_i$, and set
\[
\lambda_i=(\kappa_i\times\mathrm{id})_\#(\chi_i\tau).
\]
Each $\lambda_i$ is a finite $t$-Frostman measure on $\R^n$, with uniformly controlled mass and Frostman constant. The coordinates and cutoffs are fixed independently of the translation support of $\tau$.

Fix an integer $K>n/2$. Applying Lemma~\ref{lem:ball} in dimension $n$ to the ball $|\zeta|\le R$ and the annuli $2^jR\le|\zeta|\le2^{j+1}R$ gives
\begin{align}
\int_{\R^n}(R^2+|\zeta|^2)^{-K}|\wh\lambda_i(\zeta)|^2\,d\zeta
&\lesssim R^{n-t-2K}
\left(1+\sum_{j\ge0}2^{j(n-t-2K)}\right)\notag\\
&\lesssim R^{n-t-2K}.
\label{eq:scaled-frostman}
\end{align}
The series converges also at $t=0$, by the choice of $K$.

For the fixed test function $f$, set, in the $i$th local coordinate system,
\[
G_i(u,z)=\psi_i(\kappa_i^{-1}(u))
T_\mu f(\kappa_i^{-1}(u),z),
\]
extended by zero outside $V_i$ in the first variable. This function is smooth and compactly supported in $u$, and rapidly decreasing in $z$. Equivalence of local Sobolev norms and Proposition~\ref{prop:operator} give, for integers $0\le j\le K$,
\[
\|G_i\|_{H^j(\R^n)}
\lesssim \|T_\mu f\|_{H^j(\Isom(\R^d))}
\lesssim R^{j-\beta/2}\|f\|_2.
\]
Expanding the polynomial weight therefore yields
\begin{align}
\int_{\R^n}(R^2+|\zeta|^2)^K|\wh G_i(\zeta)|^2\,d\zeta
&\lesssim\sum_{j=0}^K R^{2(K-j)}\|G_i\|_{H^j(\R^n)}^2\notag\\
&\lesssim R^{2K-\beta}\|f\|_2^2.
\label{eq:scaled-operator}
\end{align}
Since $\psi_i=1$ on $\spt\chi_i$, Fourier inversion and Cauchy--Schwarz with the reciprocal weights in \eqref{eq:scaled-frostman} and \eqref{eq:scaled-operator} imply
\begin{align*}
\left|\int f\,d\nu\right|
&=\left|\sum_i\int G_i\,d\lambda_i\right|\\
&\le\sum_i
\left(\int(R^2+|\zeta|^2)^{-K}|\wh\lambda_i(-\zeta)|^2\,d\zeta\right)^{1/2}
\left(\int(R^2+|\zeta|^2)^K|\wh G_i(\zeta)|^2\,d\zeta\right)^{1/2}\\
&\lesssim R^{(n-t-\beta)/2}\|f\|_2.
\end{align*}
This completes
the proof.
\end{proof}
\endgroup

\subsubsection{Proof of the incidence estimate}
\begingroup
\begin{proof}[Proof of the upper bound in Theorem~\ref{thm:incidence}]

Let $\nu$ be the image of $\mu\times\tau$ under
$(x,g,z)\mapsto gx+z$. To apply the Fourier estimate in
Proposition~\ref{prop:motion-shell}, choose a nonnegative Schwartz
function $\phi$ on $\R^d$ which is at least one on the unit ball.
Majorizing the incidence indicator by $\phi$ and using Fourier inversion
gives
\begin{align*}
\mathcal I_\delta(\mu,\tau,\rho)
&\le\iint\phi((w-y)/\delta)\,d\nu(w)d\rho(y)\\
&\le\delta^d\int|\wh\nu(\xi)||\wh\rho(\xi)||\wh\phi(\delta\xi)|\,d\xi.
\end{align*}
Lemma~\ref{lem:ball}, applied to the $b$-Frostman measure
$\rho$ with weight one, gives
\[
\int_{|\xi|\le2R}|\wh\rho(\xi)|^2\,d\xi\lesssim R^{d-b}.
\]
On the annulus $2^j\le|\xi|\le2^{j+1}$, Cauchy--Schwarz and
Proposition~\ref{prop:motion-shell} therefore bound the product of the
Fourier transforms by $C2^{j(n+d-\beta-b-t)/2}$ in $L^1$.
The contribution from $|\xi|\le1$ is bounded because both measures have
mass one. The rapid decay of $\wh\phi$ now gives, for every $D>0$,
\[
\mathcal I_\delta(\mu,\tau,\rho)\lesssim_D\delta^d
\left[1+\sum_{j\ge0}2^{j(n+d-\beta-b-t)/2}(1+\delta2^j)^{-D}\right].
\]
Set $a=(n+d-\beta-b-t)/2$ and take $D>\max\{0,a\}$.
Choose the integer $J\ge1$ such that $2^J\le\delta^{-1}<2^{J+1}$.
For $0\le j\le J$, the decay factor is comparable to one. For $j>J$,
it is bounded by $(\delta2^j)^{-D}$, and $D>a$ makes this tail a
convergent geometric series. Thus,
\[
1+\sum_{j\ge0}2^{ja}(1+\delta2^j)^{-D}
\lesssim
\begin{cases}
\delta^{-a},&a>0,\\
\log(2/\delta),&a=0,\\
1,&a<0.
\end{cases}
\]
When $a>0$, the exponent after multiplication by $\delta^d$
is $d-a=(\beta+b+t-m)/2$, since $n=m+d$. The cases $a=0$ and $a<0$
give the critical and supercritical bounds, respectively. Hence,
\eqref{eq:incidence} holds in all three regimes. This completes the proof.
\end{proof}

\begin{lemma}\label{lem:incidence-basic}
Under the hypotheses of Theorem~\ref{thm:incidence}, the bounds in \eqref{eq:incidence-basic} hold.
\end{lemma}
\begin{proof}
The $b$-Frostman condition on the target measure $\rho$ gives
the first bound directly:
\[
\mathcal I_\delta(\mu,\tau,\rho)
=\iint\rho(B(gx+z,\delta))\,d\mu(x)d\tau(g,z)
\lesssim\delta^b.
\]
For the second bound, fix $x\in\spt\mu$ and
$y\in\spt\rho$. If $|gx+z-y|\le\delta$, then $(g,z)$ is within
$\delta$ of $(g,y-gx)$. Thus, the relevant motions lie in the
$\delta$-neighbourhood of
\[
\{(g,y-gx):g\in O(d)\}.
\]
The map $g\mapsto(g,y-gx)$ has Lipschitz constant bounded
in terms of $\spt\mu$, independently of $y$. Since $O(d)$ is compact
and $m$-dimensional, the displayed set, and hence its
$\delta$-neighbourhood, is covered by $O(\delta^{-m})$ balls of radius
$C\delta$. Writing their centres as $\theta_j$ and using the
$t$-Frostman condition gives
\[
\tau\{(g,z):|gx+z-y|\le\delta\}
\le\sum_{j=1}^{\lceil C\delta^{-m}\rceil}\tau(B(\theta_j,C\delta))
\lesssim\delta^{-m}\delta^t.
\]
Integrating with respect to $\mu\times\rho$ proves the second bound.

For the third bound, fix $(g,z)$ and majorize the incidence
indicator by the Schwartz function used in the preceding proof.
Orthogonal changes of variables and polar coordinates convert
\eqref{eq:spherical} into the first estimate below; applying
Lemma~\ref{lem:ball} to $\rho$ gives the second. For $R\ge1$,
\[
\int_{R\le|\xi|\le2R}|\wh\mu(g^T\xi)|^2\,d\xi\lesssim R^{d-\beta},
\qquad
\int_{|\xi|\le2R}|\wh\rho(\xi)|^2\,d\xi\lesssim R^{d-b}.
\]
Thus, Cauchy--Schwarz bounds the incidence integral for this fixed motion by
\[
C_D\delta^d\left[1+\sum_{j\ge0}2^{j(d-(\beta+b)/2)}(1+\delta2^j)^{-D}\right]
\lesssim\delta^{(\beta+b)/2}
\]
when $\beta+b<2d$, where we choose $D>d$ and split the sum
at $2^j\asymp\delta^{-1}$ as in the preceding proof. The bound is
uniform in $(g,z)$, so integration against the probability measure
$\tau$ preserves it. The only remaining case is $\beta+b=2d$:
since $\beta,b\le d$, this forces $b=d$, and the first bound already
gives $\delta^d=\delta^{(\beta+b)/2}$ without a logarithmic loss.
Combining the three bounds proves \eqref{eq:incidence-basic}.
This completes the proof.
\end{proof}
\endgroup

\subsubsection{Sharpness in every dimension}\label{subsec:incidence-sharpness}
\begingroup
\begin{proof}[Proof of the sharpness assertion in Theorem~\ref{thm:incidence}]
We construct the source measure as the product of circle measure
and a smooth density in the remaining directions. The rotations will
preserve its vertical marginal, whose square-root singularity at a
tangent hyperplane leads to the critical logarithm. We first verify the
spherical-decay and Frostman hypotheses for the construction.

Fix $d\ge2$ and $0<\alpha<1$. Let $P=\operatorname{span}(e_1,e_d)$, $Q=P^\perp$, and $V=e_d^\perp$. Choose a nonnegative $f\in C_c^\infty(Q)$ supported in the unit ball and having integral one, and set
\[
\mu=\sigma_P\otimes f(v)\,dv,
\]
where $\sigma_P$ is normalized arc length on the unit circle in $P$. When $d=2$, the second factor is the unit mass on $Q=\{0\}$.

We first verify the spherical decay. Write $\xi_P$ and $\xi_Q$ for the orthogonal projections of $\xi$ onto $P$ and $Q$. The circle Fourier estimate and rapid decay of $\wh f$ give
\[
|\wh\mu(\xi)|^2\lesssim_D
(1+|\xi_P|)^{-1}(1+|\xi_Q|)^{-D}.
\]
When $d\ge3$, the projection of normalized surface measure on $S^{d-1}$ onto $Q$ has constant density on its unit ball $B_Q(0,1)$. Explicitly, for every integrable function $\varphi$ on $S^{d-1}$,
\[
\int_{S^{d-1}}\varphi(\omega)\,d\sigma_{d-1}(\omega)
=\frac{1}{\Leb^{d-2}(B_Q(0,1))}\int_{B_Q(0,1)}\int_{S^1\subset P}
\varphi\bigl(\sqrt{1-|v|^2}\,u+v\bigr)\,d\sigma_P(u)\,dv.
\]
On $|v|\le1/2$, we have $\sqrt{1-|v|^2}\ge\sqrt3/2$;
on $|v|>1/2$, the factor $(1+R|v|)^{-D}$ is at most $C_DR^{-D}$.
Splitting the integral accordingly, for $R\ge1$ and $D>d-1$ we obtain
\begin{align*}
\sigma(\mu)(R)
&\lesssim_D\int_{B_Q(0,1)}
\frac{(1+R|v|)^{-D}}{1+R\sqrt{1-|v|^2}}\,dv\\
&\lesssim R^{-1}\int_Q(1+R|v|)^{-D}\,dv+R^{-D}
\lesssim R^{-(d-1)}.
\end{align*}
The last integral is $O(R^{-(d-2)})$ by the substitution
$v\mapsto Rv$, and $R^{-D}\le R^{-(d-1)}$. For $d=2$, the required
bound follows directly from the circle estimate, since the second
factor of $\mu$ is a point mass. Thus, $\mu$ satisfies
\eqref{eq:spherical} with $\beta=d-1$; boundedness for $0\le R\le1$
follows from its total mass.

Let
\[
H=\{h\in O(d):he_d=e_d\}\cong O(d-1),
\]
and let $R_\vartheta$ rotate $P$ through angle $\vartheta$, acting as the identity on $Q$. With $\TT=\R/(2\pi\Z)$, the map
\[
H\times\TT\longrightarrow O(d),\qquad
(h,\vartheta)\longmapsto hR_\vartheta,
\]
is a smooth embedding. Indeed, if $g=hR_\vartheta$, then $g^Te_d=R_{-\vartheta}e_d$ determines $\vartheta$ modulo $2\pi$, after which $h=gR_{-\vartheta}$. These formulas also give a smooth inverse on the image. The pushforward $\gamma$ of Haar probability on $H$ times normalized angular measure is therefore an $(m-d+2)$-Frostman probability measure, since
\[
\dim H+1=\frac{(d-1)(d-2)}2+1=m-d+2.
\]
For $d=2$, the image is all of $O(2)$ and $\gamma$ is Haar probability.

Let $\eta$ be an \emph{$\alpha$-Ahlfors regular} probability measure on a compact $C\subset[0,1/100]$; that is,
\[
\eta(B(p,r))\sim r^\alpha,
\qquad p\in C,\quad 0<r\le1.
\]
One may take a scaled copy of the equal-weight self-similar measure generated by two similarities of contraction ratio $2^{-1/\alpha}$ with disjoint images of $[0,1]$. For $a>0$, let $\lambda_a$ be normalized Lebesgue measure on $B_V(0,a)$. Define
\begin{align*}
\tau&=\bigl((g,z,p)\mapsto(g,z+pe_d)\bigr)_\#
(\gamma\times\lambda_1\times\eta),\\
\rho&=\bigl((y,p')\mapsto y+(1+p')e_d\bigr)_\#
(\lambda_4\times\eta).
\end{align*}
The product Frostman estimates give
\[
t=(m-d+2)+(d-1)+\alpha=m+1+\alpha,
\qquad b=d-1+\alpha.
\]

To prove the claimed incidence orders, we will show that
\begin{equation}\label{eq:sharp-energy}
\mathcal I_\delta(\mu,\tau,\rho)
\sim \delta^d\iint(|p-p'|+\delta)^{-1/2}\,d\eta(p)d\eta(p').
\end{equation}
We first evaluate the energy on the right to explain why this is the
required estimate. Set
\[
F(r)=(\eta\times\eta)\{(p,p'):|p-p'|\le r\}.
\]
Integrating the Ahlfors regularity bounds over $p$ gives
$F(r)\sim r^\alpha$ for $0<r\le r_0$, with a fixed $r_0>0$.
Since $|p-p'|<1$, the identity
\[
(|p-p'|+\delta)^{-1/2}
=(1+\delta)^{-1/2}
+\frac12\int_{|p-p'|}^1(r+\delta)^{-3/2}\,dr
\]
and Tonelli's theorem give, for $0<\delta<r_0$,
\begin{align*}
\iint(|p-p'|+\delta)^{-1/2}\,d\eta(p)d\eta(p')
&=(1+\delta)^{-1/2}
+\frac12\int_0^1(r+\delta)^{-3/2}F(r)\,dr\\
&\sim 1+\int_0^{r_0}r^\alpha(r+\delta)^{-3/2}\,dr\\
&\sim 1+\delta^{\alpha-1/2}
+\int_\delta^{r_0}r^{\alpha-3/2}\,dr\\
&\sim
\begin{cases}
\delta^{\alpha-1/2},&0<\alpha<1/2,\\
\log(1/\delta),&\alpha=1/2,\\
1,&1/2<\alpha<1.
\end{cases}
\end{align*}
Here the integral over $[r_0,1]$ is bounded, and the split at $r=\delta$
separates the ranges where $r+\delta$ is comparable to $\delta$ and to
$r$, respectively. Thus, \eqref{eq:sharp-energy}, once proved, gives
\eqref{eq:incidence-sharp-family} in all three cases.

It remains to derive \eqref{eq:sharp-energy} from the geometry of
the construction. The measure $\mu$ is invariant under every
$R_\vartheta$, and every $h\in H$ fixes the last coordinate. Thus, for every $g=hR_\vartheta$ in the support of $\gamma$, the last coordinate of $gx$ has the circle's vertical marginal. The variable $s=1-(gx)_d$ has density
\[
w(s)=\frac{1}{\pi\sqrt{s(2-s)}},\qquad 0<s<2.
\]
In particular, $w(s)\sim s^{-1/2}$ near zero. The horizontal component of $gx+z$ lies in $B_V(0,3)$, since $|x|\le\sqrt2$ and $|z|\le1$. For $0<\delta\le1/100$, every relevant horizontal cross-section of the incidence ball therefore lies inside $B_V(0,4)$. Integrating against the target's horizontal density gives the exact formula
\begin{equation}\label{eq:sharp-kernel}
\mathcal I_\delta(\mu,\tau,\rho)
=4^{1-d}\iint K_\delta(p-p')\,d\eta(p)d\eta(p'),
\qquad
K_\delta(u)=\int_0^2
\bigl(\delta^2-(u-s)^2\bigr)_+^{(d-1)/2}w(s)\,ds,
\end{equation}
where  $u_+=\max\{u,0\}$.

To compare \eqref{eq:sharp-kernel} with the energy in
\eqref{eq:sharp-energy}, we estimate the kernel for both signs of $u$.
For $0\le u\le1/100$, we claim that
\[
K_\delta(u)\sim \delta^d(u+\delta)^{-1/2}.
\]
For the upper bound, the integrand vanishes unless $(u-\delta)_+\le s\le u+\delta$, so
\[
K_\delta(u)\lesssim\delta^{d-1}
\int_{(u-\delta)_+}^{u+\delta}s^{-1/2}\,ds
\lesssim\delta^d(u+\delta)^{-1/2}.
\]
For the lower bound, restrict to
$u+\delta/4\le s\le u+\delta/2$. This interval lies in $(0,2)$ under
the stated bounds on $u$ and $\delta$. On it, $|u-s|\le\delta/2$, so
the cross-section factor is at least
$(3/4)^{(d-1)/2}\delta^{d-1}$, while
$w(s)\gtrsim(u+\delta)^{-1/2}$. Integration over an interval of length
$\delta/4$ proves the lower bound.

For $u\ge\delta$, $K_\delta(-u)=0$. For $0\le u<\delta$,
\[
K_\delta(-u)\lesssim\delta^{d-1}\int_0^\delta s^{-1/2}\,ds
\lesssim\delta^d(u+\delta)^{-1/2}.
\]
The preceding bounds imply
\[
K_\delta(u)+K_\delta(-u)\sim\delta^d(|u|+\delta)^{-1/2},
\qquad |u|\le1/100.
\]
Interchanging $p$ and $p'$ preserves $\eta\times\eta$, so
\[
\iint K_\delta(p-p')\,d\eta(p)d\eta(p')
=\frac12\iint\bigl(K_\delta(p-p')+K_\delta(p'-p)\bigr)
\,d\eta(p)d\eta(p').
\]
Applying the preceding two-sided kernel bound in
\eqref{eq:sharp-kernel} therefore gives the desired energy comparison.

This proves \eqref{eq:sharp-energy} and hence, by the energy
calculation above, \eqref{eq:incidence-sharp-family}. Finally,
\[
\beta+b+t=m+2d-1+2\alpha,
\]
so $\alpha<1/2$, $\alpha=1/2$, and $\alpha>1/2$ correspond to the three regimes of \eqref{eq:incidence}. The construction therefore attains each regime in every dimension, along the stated parameter family. This completes the proof.
\end{proof}
\endgroup

\subsection{Point-line incidences}\label{subsec:point-line}
\begingroup
We recover Orponen's point-line incidence inequality
\cite[Theorem~1.7]{Orponen} from Proposition~\ref{prop:operator}.
Orponen's proof uses the classical Sobolev smoothing of the \emph{X-ray
transform}. Here we obtain the required estimate by averaging a fixed
line segment over rigid motions and integrating out the translation
along the line.

For $\vartheta\in\R$, write
\[
v_\vartheta=(\cos\vartheta,\sin\vartheta),\qquad
v_\vartheta^\perp=(-\sin\vartheta,\cos\vartheta),\qquad
\ell_{\vartheta,p}=p v_\vartheta^\perp+\R v_\vartheta.
\]
The energy of a measure on the $(\vartheta,p)$-plane is computed with
its Euclidean metric, as in \eqref{eq:riesz-energy}.

\begin{proposition}[]\label{prop:point-line}
Let $\rho$ and $\eta$ be finite positive measures supported in
$B(0,1)\subset\R^2$ and $[-\pi/2,\pi/2]\times[-1,1]$, respectively.
For every $1<t<2$ and $0<\delta\le1/50$,
\begin{equation}\label{eq:point-line}
\iint\mathbf1_{\{\dist(y,\ell_{\vartheta,p})\le\delta\}}
\,d\rho(y)d\eta(\vartheta,p)
\lesssim_t\delta\,I_{3-t}(\rho)^{1/2}I_t(\eta)^{1/2}.
\end{equation}
The implicit constant is independent of $\delta$, $\rho$, and $\eta$.
\end{proposition}

\begin{proof}
We may assume that both measures are nonzero and that
both energies are finite. Set $a=(2-t)/2\in(0,1/2)$. For $f\in C_c^\infty(B(0,2))$, define its X-ray transform by
\[
\mathcal Rf(\vartheta,p)
=\int_\R f(pv_\vartheta^\perp+sv_\vartheta)\,ds.
\]
We will establish the line-integral estimate
\begin{equation}\label{eq:line-duality}
\left|\int\mathcal Rf\,d\eta\right|
\lesssim_a\|\eta\|_{H^{-a}(\R^2)}
\|f\|_{H^{a-1/2}(\R^2)}.
\end{equation}
We first explain how it gives the incidence bound.

Choose a nonnegative $\phi\in C_c^\infty(B(0,3))$ with
$\phi\ge1$ on $B(0,2)$, and set
$\phi_\delta(x)=\delta^{-2}\phi(x/\delta)$ and
$f_\delta=\rho*\phi_\delta$. Then
$\spt f_\delta\subset B(0,1+3\delta)\subset B(0,2)$.
If $\dist(y,\ell_{\vartheta,p})\le\delta$, the segment
$\ell_{\vartheta,p}\cap B(y,2\delta)$ has length at least
$2\sqrt3\,\delta$. Consequently,
\[
\mathbf1_{\{\dist(y,\ell_{\vartheta,p})\le\delta\}}
\lesssim\delta\,\mathcal R\bigl(\phi_\delta(\cdot-y)\bigr)(\vartheta,p).
\]
Assuming \eqref{eq:line-duality}, Tonelli's theorem gives
\[
\iint\mathbf1_{\{\dist(y,\ell_{\vartheta,p})\le\delta\}}
\,d\rho(y)d\eta(\vartheta,p)
\lesssim_t\delta\,
\|\eta\|_{H^{(t-2)/2}}\|f_\delta\|_{H^{(1-t)/2}}.
\]
By \eqref{eq:riesz-energy} and
$|\wh\phi(\delta\xi)|\le\|\phi\|_1$, uniformly in $\delta$,
\[
\|\eta\|_{H^{(t-2)/2}}\lesssim_t I_t(\eta)^{1/2},
\qquad
\|f_\delta\|_{H^{(1-t)/2}}
\lesssim\|\rho\|_{H^{(1-t)/2}}
\lesssim_t I_{3-t}(\rho)^{1/2}.
\]
Since $-a=(t-2)/2$ and $a-1/2=(1-t)/2$, these estimates
would prove \eqref{eq:point-line}. It remains to establish
\eqref{eq:line-duality}. We do so for every $a\ge0$ by applying
Proposition~\ref{prop:operator} to normalized length measure on a segment
and lifting $\eta$ to a measure on the motion group.

Let $\mu=\frac1{10}\mathcal H^1|_{[-5,5]\times\{0\}}$.
Direct integration gives
\[
\wh\mu(\xi)=\frac{\sin(10\pi\xi_1)}{10\pi\xi_1},
\qquad
|\wh\mu(\xi)|\lesssim\min\{1,|\xi_1|^{-1}\}.
\]
The quotient is understood by continuity at $\xi_1=0$.
For $R\ge1$, splitting the angular integral at
$|\cos\vartheta|=R^{-1}$ yields
\[
\sigma(\mu)(R)
\lesssim\int_0^{2\pi}\min\{1,(R|\cos\vartheta|)^{-2}\}\,d\vartheta
\lesssim R^{-1}.
\]
For bounded $R$, use $|\wh\mu|\le1$. Thus,
Proposition~\ref{prop:operator} with $\beta=1$ gives
\[
\|T_\mu f\|_{H^a(\Isom(\R^2))}
\lesssim_a\|f\|_{H^{a-1/2}(\R^2)},\qquad a\ge0.
\]

Choose a nonnegative $h\in C_c^\infty((-1,1))$ with
$\int h=1$, and let $R_\vartheta$ be rotation through angle
$\vartheta$. Define $\tau$ to be the pushforward
of $\eta\otimes h(u)\,du$ under
\[
(\vartheta,p,u)\longmapsto
\bigl(R_\vartheta,pv_\vartheta^\perp+uv_\vartheta\bigr).
\]
This is a smooth coordinate change on a fixed neighborhood of the
support. Indeed, the rotation determines $\vartheta$ in a coordinate
containing $[-\pi/2,\pi/2]$, and the translation $z$ determines
$p=z\cdot v_\vartheta^\perp$ and $u=z\cdot v_\vartheta$.
Boundedness of changes of coordinates on local Sobolev spaces gives
\[
\|\tau\|_{H^{-a}(\Isom(\R^2))}
\lesssim_a\|\eta\otimes h\|_{H^{-a}(\R^3)}
\le\|h\|_2\|\eta\|_{H^{-a}(\R^2)},\qquad a\ge0.
\]
The last inequality follows directly from
\begin{align*}
\|\eta\otimes h\|_{H^{-a}(\R^3)}^2
&=\int_{\R^2}\int_\R
\frac{|\wh\eta(\zeta)|^2|\wh h(\xi)|^2}
{(1+4\pi^2|\zeta|^2+4\pi^2\xi^2)^a}\,d\xi\,d\zeta\\
&\le\|h\|_2^2\|\eta\|_{H^{-a}(\R^2)}^2.
\end{align*}

If $|u|\le1$, then $[-2,2]\subset[-5+u,5+u]$. Since the integrand
vanishes for $|s|\ge2$, we have
\[
T_\mu f(R_\vartheta,pv_\vartheta^\perp+uv_\vartheta)
=\frac1{10}\int_{-5+u}^{5+u}
f(pv_\vartheta^\perp+sv_\vartheta)\,ds
=\frac1{10}\mathcal Rf(\vartheta,p).
\]
Integrating against $\eta\otimes h(u)\,du$ and applying Sobolev
duality therefore gives
\begin{align*}
\left|\int\mathcal Rf\,d\eta\right|
&=10\left|\int T_\mu f\,d\tau\right|\\
&\lesssim_a
\|\eta\|_{H^{-a}(\R^2)}\|f\|_{H^{a-1/2}(\R^2)}.
\end{align*}

This proves \eqref{eq:line-duality}; the reduction above then gives
\eqref{eq:point-line}. This completes the proof.
\end{proof}

The finite-energy assumptions in \eqref{eq:point-line} give a bound
without an additional logarithmic factor. Frostman bounds at the
corresponding exponents need not give finite energies at those
exponents. Thus, the estimate for this lifted line family is consistent
with the necessary critical logarithm in Theorem~\ref{thm:incidence}.
\endgroup

\subsection{Motions mapping a pattern into a zero set}
\begin{corollary}[]\label{cor:null}
Let $E\subset\R^d$ be nonempty and Borel, and let ${Z}_{0}\subset\R^d$ have Lebesgue outer measure zero. Then
\[
\dimH\{(g,z):gE+z\subset {Z}_{0}\}\le n-\dimSF E.
\]

\end{corollary}

\begin{proof}
\begingroup

We may assume that $\dimSF E>0$. For each $0<\beta<\dimSF E$, we will place
the motion family in a Borel family of dimension at most $n-\beta$;
letting $\beta\uparrow\dimSF E$ will then prove the claim.

Choose a $G_\delta$ set $V=\bigcap_{j\ge1}O_j$ containing $Z_0$ with $\Leb^d(V)=0$, where each $O_j$ is open. For $0<\beta<\dimSF E$, choose a compact $K_\beta\subset E$ supporting a measure satisfying \eqref{eq:spherical} with exponent $\beta$. Then
\[
\mathcal C_\beta
=\{(g,z):gK_\beta+z\subset V\}
=\bigcap_{j\ge1}\{(g,z):gK_\beta+z\subset O_j\}.
\]
Each set in the intersection is open: a compact set contained in $O_j$
remains in $O_j$ under sufficiently small perturbations of the motion.
Thus, $\mathcal C_\beta$ is $G_\delta$. If it is nonempty, then
$\dimSF K_\beta\ge\beta$ and $\mathcal C_\beta(K_\beta)\subset V$ allow
us to apply Theorem~\ref{thm:main}(ii), giving
$\dimH\mathcal C_\beta\le n-\beta$. If it is empty, there are no motions
mapping $E$ into $Z_0$. In the remaining case, every such motion belongs
to $\mathcal C_\beta$, so its family has dimension at most $n-\beta$.
Letting $\beta\uparrow\dimSF E$ completes the proof.
\endgroup
\end{proof}

\subsection{Patterns contained in a line}

\begin{proposition}\label{prop:line}
Suppose $E\subset\R^2$ and $\Theta\subset \Isom(\R^2)$ are nonempty Borel sets, $E$ lies in an affine line, $s=\dimH E\in(0,1]$, and $t=\dimH\Theta>1$. Then
\[
\dimH\Theta(E)\ge
\min\left\{s+t-1,\frac{3s+t-1}{2},s+1\right\}.
\]
This bound may be combined with Theorem~\ref{thm:main} by taking their maximum.
\end{proposition}

\begin{proof}
\begingroup
We first show that the family of lines obtained
from the rigid motions has dimension at least $t-1$.
The desired bound then follows from the Ren--Wang theorem.

By applying a fixed rigid motion, we may assume that
$E\subset\R\times\{0\}$. This does not change the union
or the Hausdorff dimension of the motion family. Let 
\[
\mathscr L=\{g(\R\times\{0\})+z:(g,z)\in\Theta\}
\]
On each orientation component, use the coordinates
\[
ge_1=v_\vartheta,\qquad z=pv_\vartheta^\perp+uv_\vartheta,
\]
where $v_\vartheta$ and $v_\vartheta^\perp$ are as above. On bounded regions where the angular intervals have length less
than $\pi$, the map $(g,z)\mapsto(\vartheta,p,u)$ is bi-Lipschitz,
and $(\vartheta,p)$ determines the image line
$\ell_{\vartheta,p}$.
We cover $\Theta$ by countably many bounded sets $\Theta_j$
on which these coordinates are well defined, with $|u|\le j$.
Let $\mathscr L_j$ be the set of image-line parameters
corresponding to $\Theta_j$.
Then
\[
\Theta_j\subset\mathscr L_j\times[-j,j],\qquad
\dimH\Theta_j\le\dimH\mathscr L_j+1
\le\dimH\mathscr L+1.
\]
The first dimension inequality uses
$\dimH(A\times[-j,j])\le\dimH A+1$ for bounded $A$.
Countable stability of Hausdorff dimension now gives
\[
\dimH\mathscr L\ge t-1.
\]
A set $F\subset\R^2$ is an \emph{$(s,t_0)$-Furstenberg set} if a family of lines of Hausdorff dimension at least $t_0$ satisfies $\dimH(F\cap\ell)\ge s$ for every line $\ell$ in the family. Each line in $\mathscr L$ contains an isometric copy of $E$ in $\Theta(E)$. Thus, $\Theta(E)$ is an $(s,t-1)$-Furstenberg set, with
$0<s\le1$ and $0<t-1\le2$. Applying the Ren--Wang theorem
\cite[Theorem 1.1]{RenWang} with line-family parameter $t-1$ gives
$\min\{s+t-1,(3s+t-1)/2,s+1\}$ as the lower bound. Both this estimate
and Theorem~\ref{thm:main} hold for the same union, so their maximum is
also a lower bound. This completes the proof.
\endgroup
\end{proof}
For example, when $\dimH E=4/5$ and $\dimH\Theta=2$, Proposition~\ref{prop:line} gives $\dimH\Theta(E)\ge17/10$. The bound obtained from $\beta_2(4/5)=1/2$ in Theorem~\ref{thm:main} is $3/2$.
\begin{proposition}\label{prop:line-area}
Let $E\subset\R^2$ and $\Theta\subset \Isom(\R^2)$ be nonempty Borel sets. If $E$ lies in an affine line and $\mathcal H^1(E)>0$, then
\[
\dimH\Theta>2\quad\Longrightarrow\quad\Leb^2(\Theta(E))>0.
\]
\end{proposition}

\begin{proof}
\begingroup
We choose a finite-energy measure on the image-line family and prove that almost all of its one-dimensional projections are absolutely continuous. Fubini's theorem then rules out a zero set containing a positive-length subset of every line in that family.

Note that the image-line family is analytic, since it is the continuous image
of the Borel motion family, and has dimension greater than one
by Proposition~\ref{prop:line}. Since the family of vertical lines
has dimension one, we may restrict to nonvertical lines.
Writing these lines as $y=a_0x+b_0$, we cover their parameters
by countably many bounded sets. At least one of these sets has
Hausdorff dimension greater than one.

Choose an exponent strictly between one and this dimension.
By Frostman's lemma, there is a compactly supported probability
measure $\eta$ on this parameter set with the chosen Frostman
exponent. In particular, $I_1(\eta)<\infty$.
For each $x$, let $\eta_x$ be the pushforward of $\eta$ under
$(a_0,b_0)\mapsto a_0x+b_0$.

For every $R_0>0$, a Fourier change of variables gives
\begin{align*}
\int_{-R_0}^{R_0}\int_\R|\wh\eta_x(\xi)|^2\,d\xi\,dx
&=\int_{-R_0}^{R_0}\int_\R|\wh\eta(x\xi,\xi)|^2\,d\xi\,dx\\
&\le C_{R_0}\int_{\R^2}|\wh\eta(\zeta)|^2|\zeta|^{-1}\,d\zeta\\
&=C_{R_0}' I_1(\eta)<\infty.
\end{align*}
The inequality uses the Jacobian $|\xi|$ and $|(x\xi,\xi)|\sim_{R_0}|\xi|$ for $|x|\le R_0$. Applying this estimate for positive integer $R_0$ shows that
$\wh\eta_x\in L^2(\R)$ for $\Leb^1$-almost every $x\in\R$.
Plancherel's theorem then gives an $L^2$ density for $\eta_x$.

Suppose that $\Theta(E)$ has Lebesgue measure zero,
and choose a Borel set $V\supset\Theta(E)$ with
$\Leb^2(V)=0$. By Fubini's theorem, for
$\Leb^1$-almost every $x$, the vertical section
$V_x=\{y:(x,y)\in V\}$ has Lebesgue measure zero,
so $\eta_x(V_x)=0$. Tonelli's theorem gives
\[
0=\int_\R\eta_x(V_x)\,dx
=\int\left(\int_\R\mathbf1_V(x,a_0x+b_0)\,dx\right)d\eta(a_0,b_0).
\]
Every line under consideration contains an isometric copy of $E$ in $V$. Since arc length along $y=a_0x+b_0$ is $\sqrt{1+a_0^2}\,dx$,
\[
\int_\R\mathbf1_V(x,a_0x+b_0)\,dx
\ge\frac{\mathcal H^1(E)}{\sqrt{1+a_0^2}}>0.
\]
The inner integral is therefore positive for every line supporting
$\eta$, contradicting its zero integral against the probability measure
$\eta$. Thus, $\Leb^2(\Theta(E))>0$, as required.
\endgroup
\end{proof}

\appendix
\section{An angular Fourier proof in the plane}\label{app:planar-fourier}
In this section, we give an angular Fourier proof of the planar case of Proposition~\ref{prop:motion-shell}. Compact support of the source measure gives exponential decay of angular modes beyond a constant multiple of the spatial frequency. Cauchy--Schwarz and the joint Frostman bound then give the required dyadic estimate. Summation recovers the planar packing bounds, and pairing with the target measure recovers the planar incidence bound.

\begin{proposition}\label{prop:planar-shell}
Let $\mu$ be a compactly supported probability measure on $\R^2$ satisfying \eqref{eq:spherical} for $0\le\beta\le2$, and let $\tau$ be a compactly supported $t$-Frostman probability measure on $\Isom(\R^2)$, where $0\le t\le3$. Define $\nu$ as in Proposition~\ref{prop:measure}. Then, for $R\ge1$,
\[
\int_{R\le|\xi|\le2R}|\wh\nu(\xi)|^2\,d\xi\lesssim R^{3-t-\beta}.
\]
\end{proposition}

\begin{proof}
We first consider orientation-preserving motions. Write their
parameters as $(\vartheta,z)\in\TT\times\R^2$, where
$\TT=\R/(2\pi\Z)$ and $R_\vartheta$ is rotation through angle
$\vartheta$. Choose $L$ so that $\spt\mu\subset B(0,L)$, and use the
mixed Fourier transform
\[
\wh\tau(k,\xi)=\int e^{-ik\vartheta-2\pi i z\cdot\xi}\,d\tau(\vartheta,z).
\]

For $k\in\mathbb Z$, let
\[
a_k(\xi)=\frac1{2\pi}\int_0^{2\pi}
\widehat\mu(R_\vartheta^T\xi)e^{-ik\vartheta}\,d\vartheta.
\]
The proof rests on two estimates. Compact support of $\mu$ will
provide constants $C_0,c>0$ such that
\[
\sum_{|k|>C_0R}|a_k(\xi)|\lesssim e^{-cR},
\qquad R\le|\xi|\le2R.
\]
The joint Frostman bound on $\tau$ will give
\begin{equation}\label{eq:mixed-ball}
\sum_{|k|\le C_0R}\int_{|\xi|\le2R}|\wh\tau(k,\xi)|^2\,d\xi
\lesssim R^{3-t}.
\end{equation}
We first explain why these estimates suffice. Assume both estimates
for the moment. The exponential tail bound makes the angular Fourier
series absolutely and uniformly convergent for each $\xi$ in the
annulus, so it can be integrated term by term against
$\tau$. Parseval's identity and the definition of $\nu$ then give
\begin{equation}\label{eq:planar-series}
\sum_{k\in\mathbb Z}|a_k(\xi)|^2
=\sigma(\mu)(|\xi|),\qquad
\widehat\nu(\xi)
=\sum_{k\in\mathbb Z}a_k(\xi)\widehat\tau(-k,\xi).
\end{equation}
Since $|\wh\tau|\le1$, we apply Cauchy--Schwarz
to the terms with $|k|\le C_0R$ in
\eqref{eq:planar-series} to get
\begin{equation}\label{eq:planar-cs}
|\wh\nu(\xi)|^2\le
C R^{-\beta}\sum_{|k|\le C_0R}|\wh\tau(-k,\xi)|^2+C e^{-cR}.
\end{equation}

Integrating \eqref{eq:planar-cs} over the annulus and using
\eqref{eq:mixed-ball} would give
\begin{equation}\label{eq:shell}
\int_{R\le|\xi|\le2R}|\wh\nu(\xi)|^2\,d\xi
\le C R^{3-t-\beta}+C R^2e^{-cR}.
\end{equation}
The exponential term is bounded by a constant multiple of
$R^{3-t-\beta}$ for $R\ge1$. Thus, it remains to establish the two
estimates above. We begin with the angular tail bound.

Fix $\xi\in\mathbb R^2$. We first extend
$\vartheta\mapsto\widehat\mu(R_\vartheta^T\xi)$ to complex values.
For $w\in\mathbb C$, define
\[
\widehat\mu(R_w^T\xi)
=\int_{\mathbb R^2}e^{-2\pi i x\cdot R_w^T\xi}\,d\mu(x),
\]
where $R_w$ is defined by the usual rotation matrix with
$\cos w$ and $\sin w$. Since $\mu$ is compactly supported,
we may differentiate under the integral on every compact
subset of $\mathbb C$. Hence, this function is entire and
$2\pi$-periodic.

Write $w=\vartheta+iu$ with $\vartheta,u\in\mathbb R$.
The formulas for the complex sine and cosine give
\[
\left|\operatorname{Im}(R_{\vartheta+iu}^T\xi)\right|
=|\xi|\,|\sinh u|.
\]
Since $\mu$ is a probability measure supported in
$B(0,L),$
\begin{align*}
\left|\widehat\mu(R_{\vartheta+iu}^T\xi)\right|
&\le\int e^{2\pi x\cdot
\operatorname{Im}(R_{\vartheta+iu}^T\xi)}\,d\mu(x)\\
&\le e^{2\pi L|\xi||\sinh u|}
\le e^{c_L|\xi|},\qquad |u|\le1,
\end{align*}
where $c_L=2\pi L\sinh 1$. This estimate is uniform in
$\vartheta$.

To obtain decay in $k$, we shift the integral defining $a_k(\xi)$
to a horizontal line in the complex plane. The function
\[
w\longmapsto\widehat\mu(R_w^T\xi)e^{-ikw}
\]
is entire and $2\pi$-periodic, since $k$ is an integer.
By Cauchy's theorem on the rectangle with vertices
$0,2\pi,2\pi+iu,iu$, the two vertical integrals cancel.
We obtain
\begin{align*}
a_k(\xi)
&=\frac1{2\pi}\int_0^{2\pi}
\widehat\mu(R_{\vartheta+iu}^T\xi)
e^{-ik(\vartheta+iu)}\,d\vartheta\\
&=\frac{e^{ku}}{2\pi}\int_0^{2\pi}
\widehat\mu(R_{\vartheta+iu}^T\xi)
e^{-ik\vartheta}\,d\vartheta.
\end{align*}

Choose $u=-1$ if $k>0$ and $u=1$ if $k<0$.
In both cases, $e^{ku}=e^{-|k|}$.
Using the preceding bound gives
\[
|a_k(\xi)|\le e^{c_L|\xi|-|k|},
\qquad k\in\mathbb Z.
\]
For $k=0$, this follows directly from the definition.
In particular, if $|k|\ge2c_L|\xi|$, then
\[
|a_k(\xi)|\le e^{-|k|/2}.
\]
Thus, the Fourier coefficients decay exponentially
when $|k|$ is sufficiently large compared with $|\xi|$.

For $R\le|\xi|\le2R$, choose $C_0>2c_L+1$ and sum the
coefficient bound as a geometric series:
\begin{align*}
\sum_{|k|>C_0R}|a_k(\xi)|
&\le2e^{2c_LR}\sum_{j=\lfloor C_0R\rfloor+1}^{\infty}e^{-j}\\
&=\frac{2e^{2c_LR-(\lfloor C_0R\rfloor+1)}}{1-e^{-1}}\\
&\le\frac{2}{1-e^{-1}}e^{-(C_0-2c_L)R}.
\end{align*}
Since $C_0-2c_L>0$, this proves the required exponential tail bound
uniformly on the annulus.

For each fixed $\xi$, the angular Fourier series
\[
\widehat\mu(R_\vartheta^T\xi)
=\sum_{k\in\mathbb Z}a_k(\xi)e^{ik\vartheta}
\]
converges absolutely and uniformly in $\vartheta$.
Indeed, the coefficients are summable, and the sum
agrees with the original continuous function by
uniqueness of Fourier coefficients. This also justifies the term-by-term integration used in
\eqref{eq:planar-series}, even when $\tau$ is singular.

\begingroup
It remains to prove \eqref{eq:mixed-ball}. We dominate the frequency
cutoff by a Gaussian so that its inverse Fourier transform is a
nonnegative kernel localized at scale $R^{-1}$.
For $|k|\le C_0R$ and $|\xi|\le2R$, we have
\[
1\le e^{C_0^2+4}e^{-(k^2+|\xi|^2)/R^2}.
\]
Fourier inversion and Poisson summation give the identity
\begin{align*}
&\sum_{k\in\Z}\int_{\R^2}e^{-(k^2+|\xi|^2)/R^2}|\wh\tau(k,\xi)|^2\,d\xi\\
&\qquad=\iint H_R(\vartheta-\vartheta',z-z')\,d\tau(\vartheta,z)d\tau(\vartheta',z'),
\end{align*}
where
\[
H_R(\varphi,w)=\pi^{3/2}R^3\sum_{j\in\Z}
\exp\left(-\frac{R^2(\varphi-2\pi j)^2}{4}-\pi^2R^2|w|^2\right).
\]
The kernel is nonnegative. Since $R\ge1$, the
periodic Gaussian satisfies
\[
\sum_{j\in\Z}
e^{-R^2(\varphi-2\pi j)^2/4}
\lesssim
e^{-cR^2\dist_{\TT}(\varphi,0)^2}.
\]
Indeed, this follows by factoring out
$e^{-R^2\dist_{\TT}(\varphi,0)^2/8}$
and observing that the remaining sum is
uniformly bounded. Consequently,
\[
H_R(\varphi,w)
\lesssim
R^3 e^{-cR^2(\dist_{\TT}(\varphi,0)^2+|w|^2)}.
\]
Using the rapid decay of the Gaussian,
we obtain, for every $D>0$,
\[
H_R(\varphi,w)\lesssim_D R^3
\bigl(1+R\dist_{\TT}(\varphi,0)+R|w|\bigr)^{-D}.
\]
Choose $D>t$ so that the dyadic annular sum converges. The Frostman
bound extends to radii greater than one after enlarging its constant,
since $\tau$ has mass one. Splitting the parameter space into a ball of
radius $R^{-1}$ and annuli of radii comparable to $2^jR^{-1}$ gives
\begin{align*}
&\sup_{(\vartheta,z)}\int H_R(\vartheta-\vartheta',z-z')\,d\tau(\vartheta',z')\\
&\qquad\lesssim R^{3-t}\left(1+\sum_{j\ge0}2^{-j(D-t)}\right)
\lesssim R^{3-t}.
\end{align*}
This proves \eqref{eq:mixed-ball}, using only the joint Frostman bound for $\tau$.
\endgroup

The two estimates have now been proved, so the initial reduction
establishes the desired bound for orientation-preserving motions.
For a general $\tau$, restrict it to the two components of $O(2)$.
Each restriction has mass at most one and satisfies the same Frostman
bound; the preceding proof applies without normalizing its mass. Fix a reflection $J$ of $\R^2$.
Every orientation-reversing orthogonal transformation
can be written as $R_\vartheta J$, where
$R_\vartheta$ is a rotation. Thus,
\[
R_\vartheta Jx+z=R_\vartheta(Jx)+z.
\]
Consequently, the contribution from orientation-reversing
motions can be treated by replacing $\mu$ with
the reflected measure $J_\#\mu$ and applying the
preceding argument to rotations. This causes no change in the estimates, since
reflection preserves spherical Fourier averages:
\[
\sigma(J_\#\mu)(R)=\sigma(\mu)(R).
\]
The Frostman condition on the motion measure is
also preserved under this change of variables.
We therefore obtain the same estimate for each
of the two components of $O(2)$.
The full image measure is the sum of the two contributions. Its
squared Fourier transform is at most twice the sum of their squared
Fourier transforms. Integrating this inequality proves the claimed
estimate for $\nu$.
This completes the proof.
\end{proof}

To recover the packing conclusions, note that $\wh\nu$ is bounded
near the origin because $\nu$ is a finite measure. If $\beta+t>3$, the
exponent $3-t-\beta$ is negative, so summing \eqref{eq:shell} over
dyadic $R$ gives $\wh\nu\in L^2$ and hence an $L^2$ density for $\nu$.
More generally, multiplying \eqref{eq:shell} by $R^{u-2}$ gives a
convergent dyadic sum whenever $u+1-t-\beta<0$. Thus,
\[
\int_{\R^2}|\wh\nu(\xi)|^2
(1+|\xi|)^{u-2}\,d\xi<\infty,
\qquad 0<u<2,\quad u<\beta+t-1.
\]
Consequently,
\[
\dimH(\spt\nu)\ge u.
\]

\section{A finite-field operator proof}\label{app:finite-field}
Throughout this appendix, $q$ is an odd prime power,
$V=\mathbb F_q^d$, $d\ge2$, and $G=O(d,\mathbb F_q)$ is the orthogonal
group of $Q(x)=x_1^2+\cdots+x_d^2$. For nonempty sets $E\subset V$ and
$\Theta\subset G\times V$, write
\[
\Theta(E)=\{gx+z:x\in E,\ (g,z)\in\Theta\}.
\]
All sums and $L^2$ norms in this appendix use counting measure. Thus,
\[
\|f\|_2^2=\sum_{y\in V}|f(y)|^2,
\qquad
\|F\|_2^2=\sum_{g\in G}\sum_{z\in V}|F(g,z)|^2.
\]
Fix a nontrivial additive character $\chi$ of $\mathbb F_q$. We use the
unnormalized Fourier transform
\[
\wh f(\xi)=\sum_{y\in V}f(y)\chi(-y\cdot\xi),
\qquad
\sum_{y\in V}|f(y)|^2=q^{-d}\sum_{\xi\in V}|\wh f(\xi)|^2.
\]

\subsection{The operator estimate}
We define the following operator
\begin{equation}\label{eq:ff-definitions}
T_Ef(g,z)=\sum_{x\in E}f(gx+z).
\end{equation}
It sums $f$ over the copy $gE+z$. The relevant Fourier quantity is
\begin{equation}\label{eq:ff-spherical-constant}
K_E=\max_{\xi\in V\setminus\{0\}}
\sum_{g\in G}|\wh{\mathbf1_E}(g^{-1}\xi)|^2.
\end{equation}
Thus, $K_E$ measures the largest Fourier square sum over the rotations of
a nonzero frequency. The zero frequency will account for the average
value of the image multiplicity and is treated separately.

\begin{proposition}\label{prop:ff-operator}
For every function $f:V\to\mathbb C$ with $\sum_{y\in V}f(y)=0$,
\begin{equation}\label{eq:ff-operator}
\sum_{g\in G}\sum_{z\in V}|T_Ef(g,z)|^2
\le K_E\sum_{y\in V}|f(y)|^2.
\end{equation}
The constant $K_E^{1/2}$ is optimal among all functions $f$
satisfying $\sum_{y\in V}f(y)=0$.
\end{proposition}

\begin{proof}
We take the Fourier transform in the translation variable to express
$\|T_Ef\|_2^2$ in terms of the sums defining $K_E$. For each $g\in G$,
the change of variables $y=gx+z$ gives
\begin{align*}
\mathcal F_z(T_Ef)(g,\xi)
&=\sum_{x\in E}\sum_{z\in V}f(gx+z)\chi(-z\cdot\xi)\\
&=\sum_{x\in E}\chi(gx\cdot\xi)
\sum_{y\in V}f(y)\chi(-y\cdot\xi)\\
&=\wh f(\xi)\wh{\mathbf1_E}(-g^{-1}\xi).
\end{align*}
Here we used $g^T=g^{-1}$. Since $\mathbf1_E$ is real-valued,
$\wh{\mathbf1_E}(-\xi)=\overline{\wh{\mathbf1_E}(\xi)}$.
Thus, Plancherel's identity in the translation variable, followed by
summation over $g$, yields

\begin{equation}\label{eq:ff-plancherel}
\sum_{g\in G}\sum_{z\in V}|T_Ef(g,z)|^2
=q^{-d}\sum_{\xi\in V}|\wh f(\xi)|^2
\sum_{g\in G}|\wh{\mathbf1_E}(g^{-1}\xi)|^2.
\end{equation}
The mean-zero assumption gives $\wh f(0)=0$. We may therefore bound the
inner sum by $K_E$ at every remaining frequency. A second application of
Plancherel gives
\[
\|T_Ef\|_2^2
\le K_Eq^{-d}\sum_{\xi\ne0}|\wh f(\xi)|^2
=K_E\|f\|_2^2,
\]
which proves \eqref{eq:ff-operator}.

To prove that the norm is optimal, choose $\xi_0\ne0$ attaining the maximum
in \eqref{eq:ff-spherical-constant} and set $f(y)=\chi(y\cdot\xi_0)$.
Character orthogonality gives
$\wh f(\xi)=q^d\mathbf1_{\{\xi_0\}}(\xi)$, so $f$ has mean zero.
Substituting into \eqref{eq:ff-plancherel} gives
$\|T_Ef\|_2^2=K_E\|f\|_2^2$. Since $f\ne0$, this completes the proof.
\end{proof}

\subsection{Duality and the packing bound}

We now introduce the family $\Theta$. Define its image multiplicity by
\begin{equation}\label{eq:ff-multiplicity}
\lambda_\Theta(y)
=\bigl|\{(x,(g,z))\in E\times\Theta:gx+z=y\}\bigr|.
\end{equation}
Each motion is injective, so $\lambda_\Theta(y)$ is the number of copies,
indexed by $\Theta$, that contain $y$. Double counting gives
\begin{equation}\label{eq:ff-pairing}
\sum_{y\in V}f(y)\lambda_\Theta(y)
=\sum_{(g,z)\in\Theta}T_Ef(g,z).
\end{equation}
This is the finite-field version of pairing the Euclidean averaging
operator with a measure on the joint parameter family. In particular,
\[
\operatorname{supp}\lambda_\Theta=\Theta(E),
\qquad
\sum_{y\in V}\lambda_\Theta(y)=|E||\Theta|.
\]
So Cauchy--Schwarz on $\Theta(E)$ therefore gives
\begin{equation}\label{eq:ff-support}
|\Theta(E)|\ge
\frac{|E|^2|\Theta|^2}{\sum_{y\in V}\lambda_\Theta(y)^2}.
\end{equation}
The role of duality is to transfer Proposition~\ref{prop:ff-operator}
to the $L^2$ bound in the following theorem.

\begin{theorem}\label{thm:ff-packing}
For all nonempty $E\subset V$ and $\Theta\subset G\times V$,
\begin{equation}\label{eq:ff-variance}
\sum_{y\in V}\lambda_\Theta(y)^2
\le \frac{|E|^2|\Theta|^2}{q^d}+K_E|\Theta|.
\end{equation}
Consequently,
\begin{equation}\label{eq:ff-packing}
|\Theta(E)|\ge
\frac{|E|^2|\Theta|}{|E|^2|\Theta|/q^d+K_E}.
\end{equation}
If $K_E>0$, then
\begin{equation}\label{eq:ff-packing-min}
|\Theta(E)|\ge\frac12
\min\left\{q^d,\frac{|E|^2|\Theta|}{K_E}\right\}.
\end{equation}
If $K_E=0$, then $\Theta(E)=V$.
\end{theorem}

\begin{proof}
We first prove the estimate below
\begin{equation}\label{eq:ff-centered}
\sum_y\left|\lambda_\Theta(y)-\frac{|E||\Theta|}{q^d}\right|^2
\le K_E|\Theta|.
\end{equation}
Indeed, since $\sum_y\lambda_\Theta(y)=|E||\Theta|$, the constant
and mean-zero parts are orthogonal, so
\[
\sum_y\lambda_\Theta(y)^2
=\frac{|E|^2|\Theta|^2}{q^d}
+\sum_y\left|\lambda_\Theta(y)-\frac{|E||\Theta|}{q^d}\right|^2.
\]
Thus, \eqref{eq:ff-centered} implies \eqref{eq:ff-variance};
substitution into \eqref{eq:ff-support} then gives \eqref{eq:ff-packing}.

For any mean-zero function $f$, \eqref{eq:ff-pairing}, Cauchy--Schwarz,
and Proposition~\ref{prop:ff-operator} give
\begin{align}
\left|\sum_y f(y)\lambda_\Theta(y)\right|
&\le |\Theta|^{1/2}
\left(\sum_{g,z}|T_Ef(g,z)|^2\right)^{1/2}\notag\\
&\le (K_E|\Theta|)^{1/2}\|f\|_2.
\label{eq:ff-duality}
\end{align}
Apply this bound with the real-valued mean-zero function
$f=\lambda_\Theta-|E||\Theta|/q^d$. For this choice,
$\sum_y f(y)\lambda_\Theta(y)=\|f\|_2^2$, and hence
\[
\|f\|_2^2\le (K_E|\Theta|)^{1/2}\|f\|_2.
\]
If $\|f\|_2=0$, \eqref{eq:ff-centered} holds immediately. Otherwise,
dividing by $\|f\|_2$ and squaring proves it. This also proves
\eqref{eq:ff-variance} and \eqref{eq:ff-packing}, by the preceding reduction.

If $K_E>0$, the denominator in \eqref{eq:ff-packing} is at most
$2\max\{|E|^2|\Theta|/q^d,K_E\}$, which gives
\eqref{eq:ff-packing-min}. If $K_E=0$, \eqref{eq:ff-centered} shows that
$\lambda_\Theta(y)=|E||\Theta|/q^d>0$ at every $y\in V$.
Thus, $\Theta(E)=V$. This completes the proof.
\end{proof}

\begin{remark}\label{rem:ff-sharp}
We can improve the bound by subtracting the average of
$\mathbf1_\Theta$. Indeed, for every mean-zero function $f$,
\[
\sum_{g,z}T_Ef(g,z)=|G||E|\sum_y f(y)=0.
\]
Thus, in \eqref{eq:ff-pairing}, we may replace
$\mathbf1_\Theta$ by
\[
\mathbf1_\Theta-\frac{|\Theta|}{|G|q^d},
\]
whose squared $L^2$ norm is
\[
|\Theta|\left(1-\frac{|\Theta|}{|G|q^d}\right).
\]
Applying Cauchy--Schwarz and Proposition~\ref{prop:ff-operator},
and then repeating the proof of Theorem~\ref{thm:ff-packing},
we obtain the improved estimate
\begin{equation}\label{eq:ff-packing-exact}
|\Theta(E)|\ge
\frac{|E|^2|\Theta|}
{|E|^2|\Theta|/q^d+
K_E\left(1-|\Theta|/(|G|q^d)\right)}.
\end{equation}
\end{remark}

\begin{remark}\label{rem:ff-weighted}
The same argument applies to arbitrary nonnegative weights $w(g,z)$;
no independence between rotations and translations is required. The
weighted multiplicity is
\[
\lambda_w(y)=\sum_{g\in G}\sum_{x\in E}w(g,y-gx).
\]
Its total mass is $|E|\sum_{g,z}w(g,z)$. Pairing with $w$ and subtracting
its average, as in Remark~\ref{rem:ff-sharp}, yields
\begin{align*}
\sum_y\left|\lambda_w(y)-\frac{|E|}{q^d}\sum_{g,z}w(g,z)\right|^2
\le K_E\left(\sum_{g,z}w(g,z)^2
-\frac{\bigl(\sum_{g,z}w(g,z)\bigr)^2}{|G|q^d}\right).
\end{align*}
Thus, the operator estimate controls the image of a joint parameter
measure directly, just as in the Euclidean proof.
\end{remark}

\subsection{The planar circle restriction estimate}

The following is the well-known circle restriction estimate for indicator
functions, stated with our unnormalized Fourier transform. It is
\cite[Lemma~4.4]{CEHIK}, also recalled in
\cite[Theorem~5.1(1)]{IMPPSS}. It follows by duality from the
$L^2$ to $L^4$ extension estimate for a circle of nonzero radius.

\begin{lemma}\label{lem:ff-circle}
Let $q$ be odd and let
$S_j=\{\xi\in\mathbb F_q^2:Q(\xi)=j\}$, where
$j\in\mathbb F_q\setminus\{0\}$. For every $E\subset\mathbb F_q^2$,
\begin{equation}\label{eq:ff-circle}
\sum_{\xi\in S_j}|\wh{\mathbf1_E}(\xi)|^2
\le\sqrt3\,q\,|E|^{3/2}.
\end{equation}
\end{lemma}

\begin{corollary}\label{cor:ff-plane}
Suppose that $q\equiv3\pmod4$. For all nonempty
$E\subset\mathbb F_q^2$ and
$\Theta\subset O(2,\mathbb F_q)\times\mathbb F_q^2$,
\begin{equation}\label{eq:ff-planar-packing}
|\Theta(E)|\gtrsim
\min\left\{q^2,\frac{|E|^{1/2}|\Theta|}{q}\right\},
\end{equation}
with an absolute implicit constant.
\end{corollary}

\begin{proof}
We will bound $K_E$ by $2\sqrt3\,q\,|E|^{3/2}$.
Theorem~\ref{thm:ff-packing} will then give
\eqref{eq:ff-planar-packing}; its case $K_E=0$ already gives
$\Theta(E)=\mathbb F_q^2$. To obtain the bound, we identify the
nonzero orthogonal orbits and their multiplicities.

Since $q\equiv3\pmod4$, $-1$ is not a square in $\mathbb F_q$.
Hence, $Q(\xi)\ne0$ whenever $\xi\ne0$.
Identify $\mathbb F_q^2$ with $\mathbb F_{q^2}$ by writing
$(x_1,x_2)$ as $x_1+ix_2$, where $i^2=-1$.
Then $Q$ is the norm map, and each set $S_j$, $j\ne0$,
contains exactly $q+1$ elements.
Moreover, the norm-one elements act transitively on $S_j$,
since $Q(\xi)=Q(\eta)=j\ne0$ implies that
$\eta/\xi$ has norm one.

We next determine the size of the orthogonal group.
Let $g\in O(2,\mathbb F_q)$. Since $Q(g(1))=1$,
multiplication by $g(1)^{-1}$ preserves $Q$.
Composing this multiplication with $g$, we obtain an
orthogonal map that fixes $1$. Such a map must preserve
the line $\mathbb F_q i$, which is perpendicular to $1$.
Since it also preserves $Q$, it sends $i$ to either
$i$ or $-i$. Thus, it is either the identity or
conjugation.

It follows that every orthogonal map is either
multiplication by a norm-one element or the composition
of such a multiplication with conjugation. Since there
are $q+1$ norm-one elements, we have
\[
|O(2,\mathbb F_q)|=2(q+1).
\]
As shown above, the orbit of every nonzero vector
$\xi$ is the set $S_{Q(\xi)}$, which contains $q+1$
elements. Therefore, each point of this orbit is
obtained from exactly two elements of the orthogonal
group. Consequently,
\[
K_E
=2\max_{j\ne0}\sum_{\xi\in S_j}
|\wh{\mathbf1_E}(\xi)|^2
\le 2\sqrt3\,q\,|E|^{3/2},
\]
where the last inequality follows from
Lemma~\ref{lem:ff-circle}.

If $K_E>0$, substituting this bound into
\eqref{eq:ff-packing-min} gives
\[
|\Theta(E)|\ge\frac12
\min\left\{
q^2,\frac{|E|^{1/2}|\Theta|}{2\sqrt3\,q}
\right\}.
\]
If $K_E=0$, Theorem~\ref{thm:ff-packing} gives
$\Theta(E)=\mathbb F_q^2$.
Thus, in both cases, we obtain
\eqref{eq:ff-planar-packing} with an absolute constant.
\end{proof}

\end{document}